\documentclass[12pt,a4paper,leqno]{article}

\usepackage{amsmath,amscd,amsfonts,amsthm,amssymb,dsfont,bbm,manfnt}
\usepackage{hyperref,url,color,appendix,eucal,comment,tabu,float}
\usepackage[usenames,dvipsnames]{xcolor}
\usepackage[nottoc,numbib]{tocbibind}
\usepackage{multicol}
\usepackage[OT2,T1]{fontenc}
\usepackage[utf8]{inputenc}
\usepackage[all]{xy}
\usepackage{lscape,rotating} 
\usepackage{graphicx}
\usepackage[nolabel]{showlabels}

\def\got{\mathfrak}

\newenvironment{pf}
{\medskip\noindent {\it Proof.  }}
{\hfill\nobreak $\Box$ \par\bigbreak}

\newcommand{\isomo}{\overset{\sim}{\rightarrow}}

\newcommand{\lisomo}{\overset{\sim}{\longrightarrow}}

\newcommand{\GL}{\mathrm{GL}}

\newcommand{\ps}{\par \smallskip}

\newcommand{\Z}{\mathbb{Z}}

\newcommand{\Q}{\mathbb{Q}}

\newcommand{\R}{\mathbb{R}}    
\newcommand{\C}{\mathbb{C}}

\newcommand{\sff}{\sffamily\selectfont}

\newtheorem*{thm*}{Theorem}
\newtheorem{lemma}[subsection]{Lemma}
\newtheorem{remark}[subsection]{Remark}
\newtheorem{cor}[subsection]{Corollary}

\newtheorem{prop}[subsection]{Proposition}
\newtheorem{example}[subsection]{Example}

\def\thesubsection{\arabic{section}.\arabic{subsection}} 

\usepackage{titlesec}

\titleformat{\section}{\bf \normalsize}{\arabic{section}.}{1 em}{}

\titleformat{\subsection}
{}
{\bf \arabic{section}.\arabic{subsection}.}
{0.5 em}
{\normalsize}

\titleformat{\subsubsection}[runin]
{\small \bf}
{}
{}
{}

\hypersetup{colorlinks=true,linkcolor=black,citecolor=black,urlcolor=Periwinkle,pdfborderstyle={/S/U/W 1}}

\makeatletter

\@addtoreset{equation}{section}
\makeatother

\title{The characteristic masses of Niemeier lattices}

\author{Ga\"etan Chenevier\thanks{Universit\'e Paris-Saclay, CNRS,  Laboratoire de math\'ematiques d'Orsay, 91405, Orsay, France. During this work, the author has been supported by the C.N.R.S. and by the projects ANR-14-CE25-0002-01 PERCOLATOR and ANR-19-CE40-0015-02 COLOSS.}}

\begin{document}
\maketitle

\begin{abstract} Let $L$ be an integral lattice in the Euclidean space $\R^n$ and $W$ an irreducible representation of the orthogonal group of $\R^n$. We give an implemented algorithm computing the dimension of the subspace of invariants in $W$ under the isometry group ${\rm O}(L)$ of $L$. A key step is the determination of the number of elements in ${\rm O}(L)$ having any given characteristic polynomial, a datum that we call the {\it characteristic masses} of $L$. As an application, we determine the characteristic masses of all the Niemeier lattices, and more generally of any even lattice of determinant $\leq 2$ in dimension $n \leq 25$. \par
For Niemeier lattices, as a verification, we provide an alternative (human) computation of the characteristic masses. The main ingredient is the determination, for each Niemeier lattice $L$ with non-empty root system $R$, of the ${\rm G}(R)$-conjugacy classes of the elements of the "umbral" subgroup ${\rm O}(L)/{\rm W}(R)$ of ${\rm G}(R)$, where ${\rm G}(R)$ is the automorphism group of the Dynkin diagram of $R$, and ${\rm W}(R)$ its Weyl group. \par
	These results have consequences for the study of the spaces of automorphic forms of the definite orthogonal groups in $n$ variables over $\Q$. As an example, we provide concrete dimension formulas in the level $1$ case, as a function of the weight $W$, up to dimension $n=25$.
\end{abstract}

\section{Introduction}

\subsection{A motivation: dimension of spaces of level $1$ automorphic forms for ${\rm O}_{n}$}\label{motivintro}
Let $n$ be an integer $ \equiv 0 \bmod 8$ and let $\mathcal{L}_n$ be the set of all {\it even unimodular} lattices in the standard Euclidean space $\R^n$. A standard example of an element of $\mathcal{L}_n$ is the lattice ${\rm E}_n={\rm D}_n + \Z \frac{1}{2}(1,1,\dots,1)$, where ${\rm D}_n$ denotes the subgroup of elements $(x_i)$ in $\Z^n$ with $\sum_i x_i \equiv 0 \bmod 2$. The orthogonal group ${\rm O}(\R^n)$ naturally acts on $\mathcal{L}_n$, with finitely many orbits, and we set 
\begin{equation} \label{defxn} {\rm X}_n\, \overset{\rm def}{=}\, {\rm O}(\R^n) \backslash \mathcal{L}_n.\end{equation} Representatives of this set ${\rm X}_n$ have been determined so far for $n\leq 24$ only: we have ${\rm X}_8=\{{\rm E}_8\}$ (Mordell), ${\rm X}_{16}=\{ {\rm E}_8 \oplus {\rm E}_8, {\rm E}_{16}\}$ (Witt) and $|{\rm X}_{24}|=24$ (Niemeier): see \cite{kneser16,Ni,venkov,conwaysloane}. The elements of $\mathcal{L}_{24}$, to which we shall refer as the {\it Niemeier lattices}, will play a major role in this paper. \ps\ps 
Similarly, for $n \equiv \pm 1 \bmod 8$ we define $\mathcal{L}_n$ as the set of all even lattices with covolume $\sqrt{2}$ in $\R^n$, as well as ${\rm X}_n$ by the same Formula \eqref{defxn}. In this case, representatives of ${\rm X}_n$ are known up to $n=121$, this last (and most complicated) case being due to Borcherds \cite{borcherds}, and we have 
$$|{\rm X}_1|=|{\rm X}_7|=|{\rm X}_9|=1, \, \, |{\rm X}_{15}|=2,\, \, |{\rm X}_{17}|=4, |{\rm X}_{23}|=32 \hspace{.3cm} \text{and} \hspace{.3cm}  |{\rm X}_{25}|=121.$$ 
\ps\ps

For any $n \equiv -1,0,1 \bmod 8$, and any complex, finite dimensional, continuous, linear representation $W$ of ${\rm O}(\R^n)$, we consider the complex vector space of $W$-valued ${\rm O}(\R^n)$-equivariant functions on $\mathcal{L}_n$: 
{\small
\begin{equation} \label{defmwo}{\rm M}_W({\rm O}_n)  = \{ f : \mathcal{L}_n \longrightarrow W\, \, \, |\,\, \, f(g L)=g f(L)\, \, \,\forall L \in \mathcal{L}_n,\, \, \,\forall g \in {\rm O}(\R^n)\}. \end{equation} 
}
\noindent This space has a natural interpretation as a space of {\it level $1$} and {\it weight $W$} automorphic forms for the orthogonal group scheme ${\rm O}_n$ of any element of $\mathcal{L}_n$. In particular, it has a very interesting action of the Hecke ring of ${\rm O}_n$ (see e.g. \cite[Sect. 4]{chlannes}), which is a first indication of our interest in it. \ps\ps

If $L$ is a lattice in the Euclidean space $\R^n$, we denote by ${\rm O}(L)=\{g \in {\rm O}(\R^n) \, \, |\, \, gL=L\}$ its (finite) isometry group. If $W$ is a representation of ${\rm O}(\R^n)$, we denote by $W^\Gamma=\{ w \in W \, \, |\, \, \gamma w = w, \, \, \forall \gamma \in \Gamma\} \subset W$ the subspace of invariants of the subgroup $\Gamma$ of ${\rm O}(\R^n)$. Fix representatives $L_1,\dots,L_h$ of the classes in ${\rm X}_n$. Then the map $f \mapsto (f(L_i))$ induces a $\C$-linear isomorphism 
\begin{equation} \label{red} {\rm M}_W({\rm O}_n) \isomo \prod_{i=1}^h W^{{\rm O}(L_i)}. \end{equation}
It follows that ${\rm M}_W({\rm O}_n)$ is finite dimensional. Our main aim in this work, which is of computational flavor, is to explain how to compute $\dim {\rm M}_W({\rm O}_{n})$ for all $n\leq 25$ and $W$ arbitrary.  The special cases $n = 7, 8, 9$ and $n=16$, more precisely their ${\rm SO}$-variants\footnote{We define ${\rm M}_{W}({\rm SO}_n)$ by replacing ${\rm O}(\R^n)$ with ${\rm SO}(\R^n)$ in \eqref{defmwo}, and $W$ with a representation of ${\rm SO}(\R^n)$. We have then ${\rm M}_W({\rm SO}_n) \simeq {\rm M}_{W'}({\rm O}_n)$ where $W'$ is the representation of ${\rm O}(\R^n)$ induced from $W$ \cite[\S 4.4.4]{chlannes}.
The question of computing dimensions in the ${\rm SO}$-case is thus a special case of the same question in the ${\rm O}$-case (the one considered here).}, had been respectively previously considered in \cite[Ch. 2]{chrenard2} and in \cite[Ch. IX Prop. 5.13]{chlannes}. In a different direction, see Appendix \ref{par:asymptotic} for an asymptotic formula for $\dim {\rm M}_W({\rm O}_n)$ (for any $n$). \ps

	Our main motivation for these computations is the relation between the spaces ${\rm M}_W({\rm O}_n)$ and geometric $\ell$-adic representations of ${\rm Gal}(\overline{\Q}/\Q)$ of Artin conductor $1$ (or pure motives over $\Q$ with good reduction everywhere) which follows from the general yoga and point of views of Langlands and Arthur on automorphic representations. This circle of ideas has been studied in great details in the recent works \cite{chlannes, chrenard2}, and pursued in \cite{taibisiegel, chetai}, to which we refer to for further explanations. As a start, the reader may consult the preface of \cite{chlannes}. Let us simply say here that in a forthcoming work of Ta\"ibi and the author, we shall use the results of the present paper as an ingredient to extend to higher dimensions $d$, hopefully up to $d=24$, the counting of level $1$, algebraic, essentially selfdual cuspidal automorphic representations of $\GL_d$ over $\Q$ started in the aforementioned works. One alternative motivating goal of these works is to obtain new information on the size of ${\rm X}_{31}$ and ${\rm X}_{32}$ (see e.g. \cite[Thm. IX.6.1]{chlannes} for a direct proof of the equality $|{\rm X}_{24}|=24$ not relying on any lattice computation). \ps\ps
	
\subsection{Dimension of invariants and characteristic masses}\label{par:introdiminv}
${}^{}$
	Consider now an arbitrary integral lattice $L$ in the standard Euclidean space $\R^n$ of arbitrary dimension $n$, and a finite dimensional representation $W$ of ${\rm O}(\R^n)$. Motivated by the previous paragraph, we are interested in algorithms to determine the dimension of 
the subspace $W^{{\rm O}(L)} \subset W$ of ${\rm O}(L)$-invariants in $W$. Of course, our requirement will be that these algorithms be efficient for the even lattices of determinant $\leq 2$, as in \S \ref{motivintro}.  \ps\ps
	Obviously, we may and do assume that $W$ is irreducible. It will be convenient to parameterize the isomorphism classes of irreducible complex representations of ${\rm O}(\R^n)$, following Weyl's original approach \cite{weyl}, by the {\it $n$-permissible}\footnote{This means that the first two columns of the Young diagram of the partition have at most $n$ boxes in total.} (integer) partitions: see Appendix \ref{par:on} for a brief reminder of this parameterization and its relation with the highest weight theory for ${\rm SO}(\R^n)$. This parameterization not only allows to deal with the two connected components of ${\rm O}(\R^{n})$ in a very concise way, but it is also especially relevant for the character formulas we shall use. \ps\ps
		We denote by ${\rm W}_\lambda$ an irreducible representation of ${\rm O}(\R^n)$ associated with the $n$-permissible partition $\lambda=(\lambda_1\geq \cdots \geq \lambda_n \geq 0)$. The element $-{\rm id}_n$ acts on ${\rm W}_\lambda$ by the sign $(-1)^{|\lambda|}$, with $|\lambda|=\sum_i \lambda_i$, so ${\rm W}_\lambda^{{\rm O}(L)}$ vanishes for  $|\lambda| \equiv 1 \bmod 2$. Our starting point is the trivial formula $\dim {\rm W}_\lambda^{{\rm O}(L)}\,=\, \frac{1}{|{\rm O}(L)|} \sum_{\gamma \in {\rm O}(L)} {\rm Trace}(\gamma ; {\rm W}_\lambda)$, that we rewrite as
\begin{equation}\label{formdim}\,\,\,\,\dim {\rm W_\lambda}^{{\rm O}(L)}\, = \, \sum_{P \in {\rm Car}_n} \,\,{\rm m}_{{\rm O}(L)}(P)\,\,\,{\rm Trace}(\,{\rm c}_P \,;\, {\rm W}_\lambda\,)\end{equation}
where: \ps\ps

(i) ${\rm Car}_n \subset \Z[t]$ denotes the (finite) subset of polynomials of degree $n$ which are products of cyclotomic polynomials. This subset is\footnote{Set $\zeta=e^{\frac{2i\pi}{m}}$ for $m\geq 1$. The symmetric bilinear form $(x,y) \mapsto {\rm Trace}_{\Q(\zeta)/\Q}( x \overline{y})$ on the free abelian group $L=\Z[\zeta]$ defines a inner product on $L \otimes \R$. The multiplication by $\zeta$ is an isometry preserving $L$, with characteristic polynomial the $m$-th cyclotomic polynomial.} also the set of characteristic polynomials of the elements of ${\rm O}(\R^n)$ preserving some lattice in $\R^n$. Using the irreducibility of cyclotomic polynomials in $\Q[t]$, it is straightforward to enumerate the elements of ${\rm Car}_{n}$ for small $n$ with the help of a computer: see Table \ref{tab:cardcarn} for the cardinality of ${\rm Car}_{n}$ for $n \leq 27$ (sequence A120963 on \cite{OEIS}). \ps\ps

\begin{table}[htp]
{\scriptsize \renewcommand{\arraystretch}{1.8} \medskip
\begin{center}
\begin{tabular}{c|c|c|c|c|c|c|c|c|c}
$n$ & $1$ & $2$ & $3$ & $4$ & $5$ & $6$ & $7$ & $8$ & $9$ \cr
\hline
$|{\rm Car}_{n}$| & $2$ & $6$ & $10$ & $24$ & $38$ & $78$ & $118$ & $224$ & $330$ \cr
\noalign{\hrule height 1pt}
$n$ & $10$ & $11$ & $12$ & $13$ & $14$ & $15$ & $16$ & $17$ & $18$ \cr
\hline
$|{\rm Car}_{n}|$ & $584$ & $838$ & $1420$ & $2002$ & $3258$ & $4514$ & $7134$ &  $9754$  & $15010$ \cr
\noalign{\hrule height 1pt}
$n$ & $19$ & $20$ & $21$ & $22$ & $23$ & $24$ & $25$ & $26$ & $27$ \cr
\hline
$|{\rm Car}_{n}|$  & $20266$   & $30532$  &   $40798$ & $60280$  &  $79762$ &   $115966$&  $152170$ & $217962$ & $283754$ \cr
\end{tabular} 
\end{center}
} 
\caption{{\small The cardinality of ${\rm Car}_{n}$ for $n\leq 27$.} }\label{tab:cardcarn}
\end{table}

(ii) For any finite subset $S \subset {\rm O}(\R^n)$, and any $P$ in $\R[t]$, we denote by ${\rm m}_S(P)$ the number of elements $g$ in $S$ with $\det(t \,{\rm id}_n -g)=P$, {\bf divided by $|S|$}. This is an element of $\Q_{\geq 0}$ that we call the {\it mass} of $P$ in $S$.  By definition, we have $$\sum_{P \in \R[t]} {\rm m}_S(P)=1.$$ 

(iii) For $P$ in $\R[t]$ a monic polynomial of degree $n$ whose complex roots are on the unit circle (e.g. $P \in {\rm Car}_n$), we denote by ${\rm c}_P \subset {\rm O}(\R^n)$ the unique conjugacy class whose characteristic polynomial is $P$. \ps\ps

We now discuss the problem of evaluating Formula \eqref{formdim}. The main unknown, which contains all the required information about $L$ and which does not depend on $\lambda$, is of course the collection of masses ${\rm m}_{{\rm O}(L)}(P)$ for $P$  in ${\rm Car}_n$. This collection will be called {\it the characteristic masses of $L$}, or sometimes simply\footnote{Beware not to confuse {\it the masses of $L$} in this sense with {\it the mass of the genus of $L$}, which traditionally appears in the study of the Minkowski-Siegel-Smith mass formula.} {\it the masses of $L$}, and we will go back to it later. We rather discuss first the question of evaluating, given an arbitrary polynomial $P$ as in (iii), the quantity ${\rm Trace}(\,{\rm c}_{P}\,; \,{\rm W}_\lambda)$. This question does not depend on $L$. \ps\ps

{\bf Evaluation of  ${\rm Trace}(\,{\rm c}_{P}\,; \,{\rm W}_\lambda)$.}  We will use for this the 
``determinantal'' character formula for ${\rm W}_{\lambda}$ proved by Weyl in \cite[Chap. VII \S 9]{weyl}. 
This formula applies to arbitrary elements of ${\rm O}(\R^{n})$, possibly of determinant $-1$. We found it useful to actually use the following alternative expression proved by Koike and Terada in \cite{koiketerada} in the spirit of the famous Jacobi-Trudi formula for the Schur polynomials in terms of elementary symmetric polynomials (see Appendix \ref{par:on}). Write $t^nP(1/t)=\sum_{i \in \Z} (-1)^{i}e_i t^i$  (so $e_i=0$ for $i<0$ or $i>n$). Denote by $\mu_1 \geq \mu_2 \geq \dots \geq \mu_{m}$ with $m=\lambda_{1}$ the partition which is dual to $\lambda$, and set $\delta_1=0$ and $\delta_j=1$ for $j>1$. Then we have the equality
\begin{equation}\label{koiketerrada} {\rm Trace} ( \,{\rm c}_P \,;\, {\rm W}_\lambda\,) \,=\, \det (e_{\mu_i-i+j}\,+\, \delta_{j}\,e_{\mu_i-i-j+2})_{1 \leq i,j \leq m} \end{equation} 
 This formula is clearly efficient when $m=\lambda_{1}$ is small, which suits well for instance the application to $|{\rm X}_{32}|$ mentioned in \S \ref{motivintro}, as it requires all $\lambda$'s with $\lambda_{1} \leq 4$ for $n=24$. Let us note that in this range, the use of the crude {\it degenerate Weyl character formula} as in \cite[\S 2]{chrenard2} would be impracticable as the Weyl group of ${\rm SO}(\R^{24})$ is much too big. Actually, the whole tables of invariants obtained in \cite[\S 2]{chrenard2} for the subgroup of determinant $1$ elements in the Weyl groups of type ${\bf E}_7$, ${\bf E}_8$ and ${\bf E}_8 \coprod {\bf A}_1$ (with respectively $n=7,8,9$) can be recomputed essentially instantly using rather Formula \eqref{koiketerrada}. \ps\ps

{\bf Determination of the characteristic masses of $L$}. This is the remaining and most important\footnote{It is equivalent to determine the finitely many ${\rm m}_{{\rm O}(L)}(P)$ for all $P$ in ${\rm Car}_n$, and the $\dim {\rm W}_\lambda^{{\rm O}(L)}$ for all $\lambda$, as the ${\rm Car}_n \times \Lambda$-matrix $({\rm Trace}(\,{\rm c}_{P}\,; \,{\rm W}_\lambda))_{P,\lambda}$ has rank $|{\rm Car}_n|$ for general reasons.
} unknown. In dimension $n$ as large as $24$, it is impossible in general to enumerate the elements of ${\rm O}(L)$ with a computer, hence to naively list their characteristic polynomials. For instance when $L$ is a Niemeier lattice then the size of ${\rm O}(L)$ is always at least $10^{14}$, and it is about $10^{30}$ for $L={\rm E}_{24}$. However, those groups have of course much fewer conjugacy classes. Write $${\rm Conj} \,{\rm O}(L) \,=\, \{ {\rm c}_i(L) \}_{i \in I}$$ 
the set of conjugacy classes of ${\rm O}(L)$. Assuming that we know representatives of the ${\rm c}_i(L)$, as well as each $|{\rm c}_i(L)|$, then the enumeration of the characteristic polynomials of ${\rm O}(L)$ may become straightforward. Of course, if we do not know representatives of ${\rm c}_i(L)$, but still the trace of the latter in $\R^n$ as well as the power maps on the ${\rm c}_i(L)$, this may similarly allow to determine the 
 characteristic masses of $L$.  \ps \ps

\begin{example} \label{Leechtable}{\rm (Leech lattice) Consider for instance the case where $L={\rm Leech}$ is ``the'' Leech lattice in $\R^{24}$. The group ${\rm O}({\rm Leech})$ is the Conway group ${\rm Co}_0$, also denoted $2\,.\, {\rm Co}_1$ in the $\mathbb{ATLAS}$ p. 180. The character of its natural representation on $\R^{24}$ is the character $\chi_{102}$ in the table {\it loc. cit.} This character, as well as Newton's relations and the power maps of the $\mathbb{ATLAS}$ (implemented in \texttt{GAP}), allow to compute the characteristic polynomial of each conjugacy class in ${\rm O}({\rm Leech})$, hence the characteristic masses of ${\rm Leech}$: they are gathered in Table \ref{TableLeech} (see p.~\pageref{not:tables} below for the notations). Note that despite the huge order $\simeq 8 \cdot 10^{18}$ of ${\rm O}({\rm Leech})$, this group only has $167$ conjugacy classes, and $160$ distinct characteristic polynomials. This is actually the minimum for a Niemeier lattice, and makes the table above printable. An interesting consequence of this computation is the observation 
$$\frac{1}{|{\rm O}({\rm Leech})|}\sum_{g \in {\rm O}({\rm Leech})} \det(\, t\, {\rm id}_{24} - g )  = t^{24}+t^{16}+t^{12}+t^8+1.$$
This asserts the existence of a line of ${\rm O}({\rm Leech})$-equivariant alternating $g$-multilinear form ${\rm Leech}^g \rightarrow \Z$ for each $g$ in $\{8,12,16,24\}$. We refer to \cite{chetaib} for a study of these forms and of the weight $13$ pluriharmonic Siegel theta series for ${\rm Sp}_{2g}(\Z)$ that they allow to construct. The results of this paper suggest several other intriguing constructions to study in the same spirit, for instance whenever a $1$ appears as a dimension for ${\rm M}_{{\rm W}_\lambda}({\rm O}_{24})$ in Table \ref{tab:dimO243}  (the case discussed here corresponding to $\lambda=\emptyset$, $1^8$ and $1^{12}$).
}
\end{example}
\ps\ps

\subsection{Algorithms for computing characteristic masses} \label{par:introalgo}

Let us give now a first algorithm, called {\sff Algorithm A} in the sequel, which takes as input the Gram matrix $G$ of some $\Z$-basis of $L$ and returns for each conjugacy class ${\rm c}_i(L)$ some representative and its cardinality $|{\rm c}_i(L)|$, hence in particular the characteristic masses of $L$. The idea, certainly classical in computational group theory, is to: \ps\ps
{\sff
\small
A1.~Apply the Plesken-Souvignier algorithm \cite{pleskensouvignier} to $G$ (implemented {\it e.g.} as $\texttt{qfauto}(G)$ in \texttt{PARI/GP}) to obtain a set $\mathcal{G}$ of generators of ${\rm O}(L)$, \ps\ps
A2. Choose a (small) finite subset $\mathcal{S} \subset L$ stable under ${\rm O}(L)$, generating $L \otimes \R$, and view ${\rm O}(L)$ as the subgroup of permutations of $\mathcal{S}$ generated by $\mathcal{G}$, \ps\ps
A3. Apply permutation groups algorithms implemented in \texttt{GAP} (such as \cite{hulpke}) to deduce cardinality and representatives of the conjugacy classes of ${\rm O}(L)$. \ps\ps
}
A canonical choice of $\mathcal{S}$ is the following: for any lattice $L$ set\label{def:S(L)} (inductively) $S(L)=M(L) \coprod S(L')$ where $M(L)$ is the subset of elements of $L$ with minimal nonzero length, and where $L'$ is the orthogonal of $M(L)$ in $L$. The choice $\mathcal{S}=S(L)$ has proved efficient enough for us in practice. We will say more about a \texttt{PARI/GP} implementation of the whole algorithm later, when discussing an improvement of it: see \S \ref{par:implementalgo}. \ps \ps

{\sff Algorithm A} is very efficient in small dimension. For instance, when $L$ is a root lattice of type ${\bf E}_6, {\bf E}_7$ or ${\bf E}_8$, it returns the characteristic masses of $L$ in a few seconds only.\footnote{All the computations in this paper have been made on a processor \texttt{Intel(R) Xeon(R) CPU E5-2650 v4 @ 2.20GHz} with \texttt{65 GB} of memory. Nevertheless, all the computations involving either {\sff Algorithm B}, or {\sff Algorithm A} in small dimension, 
are equally efficient on our personal computer (processor \texttt{1,8 GHz Intel Core i5} with \texttt{8 GB} of memory).} It turns out that it is still terminates for most of the even lattices of determinant $\leq 2$ and dimension $\leq 25$, with running time varying from a few minutes to a few days in dimensions $23,24$ and $25$ when terminates. For instance, in the case $L={\rm Leech}$ it allows to re-compute Table \ref{TableLeech} from scratch, without relying at all on the $\mathbb{ATLAS}$: it requires about $3$ minutes for step {\sff A1}, nothing for {\sff A2}, and $42$ minutes for {\sff A3}. On the other hand, it does not terminate for instance on our computer for lattices $L$ in $\mathcal{L}_{25}$ with root system\footnote{For $n\leq 25$, it follows from the classification of ${\rm X}_n$ recalled in \S \ref{motivintro} that two lattices in $\mathcal{L}_n$ are isometric if, and only if, they have isomorphic root systems.}
 ${\bf A}_1 \,{\bf D}_4\, 2{\bf D}_6\, {\bf D}_8$ or ${\bf A}_1 \,{\bf D}_6\, {\bf D}_8\,{\bf D}_{10}$ (memory issue). {\sff Algorithm A} is typically very slow (and memory consuming) if either $L$ has too many vectors $v$ of length $v \cdot v = G_{i,i}$ for some $i=1,\dots,n$, because of step {\sff A1}, or if ${\rm O}(L)$ has too many conjugacy classes, because of step {\sff A3}. It is also quite sensitive to the choice of Gram matrix $G$ of $L$ in step {\sff A1}. \ps


 In \S \ref{par:genalgo}, we will explain a significant improvement of {\sff Algorithm A} when $L$ has a non trivial root system. The basic idea of this {\sff Algorithm B} is to first write
 $${\rm O}(L) = {\rm W}(R) \rtimes {\rm O}(L)_\rho$$
where $R$ is the root system of $L$, ${\rm W}(R)$ its Weyl group, $\rho$ a Weyl vector of $R$ and ${\rm O}(L)_\rho$ the stabilizer of $\rho$ in ${\rm O}(L)$.
As we shall see, we may actually reduce the computation of the characteristic masses of $L$ to that of representatives $\gamma_j$, and sizes, of the conjugacy classes of the smaller group ${\rm O}(L)_\rho$, an information which can be obtained by replacing ${\rm O}(L)$ with ${\rm O}(L)_\rho$ in steps {\sff A1} and {\sff A3} of {\sff Algorithm A}. There are two ingredients for this reduction. The first is the determination, for each rank $r$ irreducible root system $R'$ of type ${\bf ADE}$, of the map ${\rm m}_S : {\rm Car}_r \rightarrow \Q$, where $S$ is any coset of ${\rm W}(R')$ in the full isometry group ${\rm O}(R')$ of the root system $R'$: see \S \ref{par:rootlat} for this step (which does not depend on $L$). The second is the determination, for each $j$, of the conjugacy class of $\gamma_j$ viewed as an element of the automorphism group of the Dynkin diagram of $R$. See \S \ref{par:genalgo} for a detailed discussion of {\sff Algorithm B} and of its implementation. \ps\ps

\begin{remark} {\rm (Generalizations)  In this paper, we use a restricted notion of root which suits well our applications to the  lattices in $\mathcal{L}_n$. A minor modification of {\sff Algorithm B} allows to consider the most general roots, namely the elements $\alpha$ of a lattice $L$ such that the orthogonal symmetry about $\alpha$ preserves $L$. In a different direction, it would be useful to extend the algorithms above to the context of hermitian or quaternionic positive definite lattices, possibly over totally real number fields, using the theory of complex or quaternionic reflection groups (see e.g. \cite{cohen1,cohen2}). That should help extending to higher ranks and weights the computations of dimension spaces of automorphic forms for definite unitary groups (hermitian or quaternionic) started in the literature (e.g. in \cite{lanskypollack, loeffler, dummigan, dembele, greenbergvoight}). } \ps\ps
\end{remark}
 
 \subsection{Main results} 

Using {\sff Algorithm B}, it only takes a few seconds to the computer to compute all the characteristic masses of each Niemeier lattices with roots, except in the case (trivial anyway) of ${\rm D}_{24}$ for which the Plesken-Souvignier algorithm needs about $2$ minutes. It is equally efficient in any dimension $\leq 25$: the characteristic masses are computed in a few seconds, except for ten lattices (in dimension $23$ or $25$) for which it requires less than $5$ minutes, and for the lattice ${\rm A}_1 \oplus {\rm Leech}$ in $\mathcal{L}_{25}$ (about $35$ minutes). We refer to the homepage \cite{homepage} for the gram matrices we used in our computations. Our main result is then the following.

\begin{thm*} 
\begin{itemize} 
\item[(i)] Assume $n\leq 25$. The characteristic masses of all $L \in \mathcal{L}_n$ are those given\footnote{They cannot be printed here: there are $53204$ polynomials $P$ with ${\rm m}_{{\rm O}(L)}(P) \neq 0$ for some $L$ in $\mathcal{L}_{24}$, that is about half $|{\rm Car}_{24}|$.
} in {\rm \cite{homepage}}.
\item[(ii)] The nonzero values of $\dim {\rm M}_\lambda({\rm O}_{24})$ for $\lambda_{1} \leq 3$ are given in Table \ref{tab:dimO243}. 
\end{itemize}
\end{thm*}

Table \ref{tab:dimO243} is deduced from assertion (i) for $n=24$ using observation (b) and Formulas \eqref{formdim} \&  \eqref{koiketerrada}. This step is very efficient: once the masses in (i) are computed, it takes only $5$ minutes about to produce this table. The format of the table is as follows. The notation $n_1^{m_1} \dots n_r^{m_r}$ for a partition $\lambda$ means that the diagram of $\lambda$ has exactly $m_i$ rows of size $n_i$ for $i=1,\dots,r$, and no other row. Set ${\rm d}_\lambda = \dim {\rm M}_{{\rm W}_\lambda}({\rm O}_{24})$ and denote by ${\rm ass}(\lambda)$ the associate of $\lambda$ (see \S \ref{par:on}). The column $\dim$ gives the integer ${\rm d}_\lambda$ in the case $\lambda = {\rm ass}(\lambda)$, and the two integers ${\rm d}_\lambda : {\rm d}_{{\rm ass}(\lambda)}$ otherwise. See \cite{homepage} for more extensive tables, including for instance all $\lambda$ with $\lambda_1=4$ and arbitray $n \leq 25$. \ps\ps

\begin{remark}\label{rem:spinnorm}  {\rm 
Let $L$ be a lattice in $\R^n$, fix $\gamma$ in ${\rm O}(L)$ and write $\det(t - \gamma)=(t-1)^a (t+1)^b Q(t)$ with $Q$ in $\Z[t]$ and $Q(-1)Q(1) \neq 0$. Assuming furthermore $n\equiv -1,0,1 \bmod 8$ and $L \in \mathcal{L}_n$ then Proposition 3.7 in \cite{chetai} 
shows\footnote{If  $\det \gamma =-1$ (so $b$ is odd) and $a=0$ (so $n$ is odd), apply 
the proposition to $-\gamma$.}
that for $a=0$ (resp. $b=0$) the integer $Q(1)$ (resp. $Q(-1)$) is a square. This constraint is in agreement with our computations.
}
\end{remark}

\subsection{A direct computation in the case of Niemeier lattices} \label{intro:niemeierverif}${}^{}$

In section \ref{par:niemeier}, we will explain an alternative (and human) computation of the characteristic masses of Niemeier lattices. By the results of \S \ref{par:rootlat}, we are left to determine, for each Niemeier lattice $L$ with non-empty root system $R$, the ${\rm G}(R)$-conjugacy classes of the elements of the subgroup ${\rm O}(L)/{\rm W}(R)$ of ${\rm G}(R)$, where ${\rm G}(R)$ is the automorphism group of the Dynkin diagram of $R$. We do so using a tedious case by case analysis. \ps

We found it useful to gather first in section \ref{sec:hn} some elementary results about the {\it hyperoctahedral} group ${\rm H}_n = \{ \pm 1\}^n \rtimes {\rm S}_n$. This group is both a typical direct summand of the ${\rm G}(R)$ above, and closely related to the Weyl groups of type ${\bf D}_n$ studied in \S \ref{par:irrmasses}. In particular, we introduce and characterize directly in \S \ref{parwr} and \S \ref{parhn} a few specific subgroups of ${\rm H}_n$ that will play a role in the analysis of Niemeier lattices in section \ref{par:niemeier}. \ps

Although more interesting (at least to us) from a mathematical point of view, it will be eventually clear that this nonautomatized method is too complicated to be used systematically: it would even require some work to attack the dimensions $23$ and $25$ along the same lines. Nevertheless, it provides an important check that the masses returned by the implementation of our algorithms are correct.\ps \bigskip

{\sc Aknowledgements:} We thank Jean Lannes and Olivier Ta\"ibi for useful discussions, the LMO for sharing the machine \texttt{pascaline},
as well as Bill Allombert for answering our questions on \texttt{PARI/GP}.

{\small
\tableofcontents
}
\bigskip

\begin{center}{\sc General notations and conventions}\end{center}  \label{notations}\ps

		In this paper, all group actions will be on the left. We denote by $|X|$ the cardinality of the set $X$. For $n\geq 1$ an integer, we denote by ${\rm S}_n$ the symmetric group on $\{1,\dots,n\}$, by ${\rm Alt}_n \subset {\rm S}_n$ the alternating subgroup, and we set $\Z/n:=\Z/n\Z$.\ps\ps
	If $V$ is an Euclidean space, we usually denote by $x \cdot y$ its inner product, with associated quadratic form ${\rm q}: V \rightarrow \R$ defined by ${\rm q}(x) = \frac{ x\cdot x}{2}$. A {\it lattice} in $V$ is a subgroup generated by a basis of $V$, or equivalently, a discrete subgroup $L$ with finite covolume, denoted ${\rm covol}\, L$. \ps \ps
	If $L$ is a lattice in the Euclidean space $V$, its {\it dual lattice} is the lattice $L^\sharp$ defined as $\{ v \in V\, \, |\, \, v \cdot x \in \Z ,\,\, \, \forall x \in L\}$. We say that $L$ is {\it integral} (resp. {\it even}) if we have $L \subset L^\sharp$ (resp. ${\rm q}(L) \subset \Z$). An even lattice is integral. If $L$ is integral, we have $({\rm covol}\, L)^2 \,=\,  |L^\sharp/L|$. This integer is also the {\it determinant} $\det L$ of the {\it Gram matrix} ${\rm Gram} (e) = (e_i \cdot e_j)_{1 \leq i,j \leq n}$ of any $\Z$-basis $e=(e_1,\dots,e_n)$ of $L$.  The orthogonal group of $L$ is the finite group ${\rm O}(L)=\{ \gamma \in {\rm O}(V), \, \, \gamma(L)=L\}$. \ps\ps
	\label{not:tables} In the tables of Appendix \ref{sec:tables}, we use the notation
$1^{a_1}\,2^{a_2}\, \dots\, m^{a_m}$ for the polynomial $\varphi_1^{a_1} \varphi_2^{a_2} \cdots \varphi_m^{a_m}$, where $\varphi_n$ is the $n$-th cyclotomic polynomial and where the symbol "$i^a$" is omitted for $a=0$, and shorten as "$i$" for $a=1$. 
	
\section{Preliminaries on the hyperoctahedral groups} \label{sec:hn}

\subsection{The hyperoctahedral group} \label{par:defhn}${}^{}$

Let $n\geq 1$ be an integer. The symmetric group ${\rm S}_n$ on the set $\{1,\dots,n\}$ acts on the elementary abelian $2$-group $\{\pm 1\}^n$ by permuting coordinates. The {\it hyperoctahedral group} on $n$ letters is defined as the semi-direct product $${\rm H}_n = \{\pm 1\}^n \rtimes {\rm S}_n.$$
Equivalently, ${\rm H}_n$ is the wreath product $\{\pm 1\}\,\, \wr\,\, S_n$. It is isomorphic to several familiar groups: the Weyl group of a root system of type ${\bf B}_n$ or ${\bf C}_n$, the subgroup of monomial matrices in ${\rm GL}_n(\Z)$, the orthogonal group of the standard unimodular lattice ${\rm I}_n$, the subgroup of the symmetric group on $\{\pm 1,\pm 2,\dots,\pm n\}$ of permutations $\sigma$ with $\sigma(-i)=-\sigma(i)$ for all $i$, {\rm etc}... \ps\ps

In this paper, we will encounter ${\rm H}_n$ first when discussing ${\rm O}({\rm D}_n)$ and again when studying automorphism groups of isotypic root systems. Certain subgroups of the hyperoctahedral groups will play a role in the study of Niemeier lattices. Here is an example of an interesting subgroup that will occur in the case $n=4$. We denote by $\pi : {\rm H}_n \rightarrow {\rm S}_n$ the canonical projection.  \ps\ps

  \begin{example}\label{exgl2z3} {\rm The group ${\rm GL}_2(\Z/3)$ acts on the $8$-elements set $(\Z/3)^2-\{0\}$ by permuting the $4$ disjoint pairs of the form $\{v,-v\}$. By the universal property of wreath products, the choice of elements $v_1,v_2,v_3,v_4$ such that $(\Z/3)^2-\{0\}=\coprod_i \{v_i,-v_i\}$ defines an embedding $\iota : {\rm GL}_2(\Z/3) \longrightarrow {\rm H}_4$ (a different choice leading to an ${\rm H}_4$-conjugate embedding). We have $\iota(-{\rm Id}_2)=-1$; the morphism $\pi \circ \iota$ has kernel $\pm {\rm Id}_2$ and induces ``the'' exceptional isomorphism ${\rm PGL}_2(\Z/3) \simeq {\rm S}_4$. The restriction of $\pi \circ \iota$ to the stabilizer of $v_i$ in ${\rm GL}_2(\Z/3)$ is an isomorphism onto the stabilizer ($\simeq {\rm S}_3$) of $i$ in $\{1,2,3,4\}$.} \end{example}\ps\ps
 
 We end this paragraph with a few notations and remarks about the basic structure of ${\rm H}_n$. We denote by $\varepsilon_i$ the element of $\{\pm 1\}^n$ whose $j^{\rm th}$-component is $1$ for $j \neq i$ and $-1$ for $j=i$. The center of ${\rm H}_n$ is generated by the element $-1= \prod_{i=1}^n \varepsilon_i$. The signature $\epsilon : {\rm S}_n \rightarrow \{\pm 1\}$, composed with the natural projection $\pi : {\rm H}_n \rightarrow {\rm S}_n$, defines a morphism ${\rm H}_n \rightarrow \{\pm 1\}$ that we will still denote by $\epsilon$. Another important morphism ${\rm s} : {\rm H}_n \rightarrow \{\pm 1\}$ is defined by 
\begin{equation}\label{defshn}{\rm s}( v \sigma ) \,\,=\,\, \prod_{i=1}^n v_i,\hspace{.3cm} \text{for\,all}\, \,\sigma \in {\rm S}_n \, \, \text{and}\,\, v=(v_i) \in \{\pm 1\}^n.\end{equation}
The product character $\epsilon \, {\rm s}$ coincides with the determinant when we view ${\rm H}_n$ as a the subgroup of monomial matrices in ${\rm GL}_n(\Z)$.

\subsection{Conjugacy classes of ${\rm H}_n$}\label{par:conjhn}${}^{}$

Let $\Sigma$ be a nonempty subset of  $\{1,\dots,n\}$. A {\it cycle} in ${\rm H}_n$ with support $\Sigma$ is an element of the form $h = v c$, where $c \in {\rm S}_n$ permutes transitively the elements of $\Sigma$ and fixes its complement, and where $v=(v_i) \in \{\pm 1\}^n$ satisfies $v_i=1$ for $i \notin \Sigma$. Such a cycle has a {\it length} ${\rm l}(h)$ defined as $|\Sigma|$, and a {\it sign}  ${\rm s}(h)$ (an element in $\{\pm 1\}$). This sign is also the $i$-th coordinate of $h^{{\rm l}(h)}$ for any $i$ in $\Sigma$, and ${\rm l}(h)$ is the order of $c$.
One easily checks that two cycles are conjugate in ${\rm H}_n$ if, and only if, they have the same length and the same sign. \ps\ps

Just as for ${\rm S}_n$, any element $h$ of ${\rm H}_n$ may be written as a product of cycles $h_i$ with disjoint supports, this decomposition being unique up to permutation of those cycles. The sum of the lengths of the cycles $h_i$ with ${\rm s}(h_i)=1$ (resp. ${\rm s}(h_i)=-1$) is an integer denoted ${\rm n}_{+}(h)$ (resp ${\rm n}_{-}(h)$); the collection of the length ${\rm l}(h_i)$ of those $h_i$ defines a integer partition of ${\rm n}_{+}(h)$ (resp. ${\rm n}_{-}(h)$) that we denote by ${\rm p}_{+}(h)$ (resp. ${\rm p}_{-}(h)$). We have ${\rm n}_{+}(h)+{\rm n}_{-}(h)=n$. The {\it type} of $h$ is defined as the couple of integer partitions $({\rm p}_{+}(h),{\rm p}_{-}(h))$. Two elements of ${\rm H}_n$ are conjugate if, and only if, they have the same type. \ps\ps

In the sequel, we will have to determine the type of all the elements of certain specific subgroups $G \subset {\rm H}_n$. For instance, when $G$ is the group $\iota({\rm GL}_2(\Z/3))$ of Example \ref{exgl2z3}, this information is given in Table \ref{tablegl2z3}, the row \texttt{size} giving the number of elements of the corresponding type {\bf divided by $|G|$}: 
\begin{table}[htp]
{\scriptsize \renewcommand{\arraystretch}{1.8} \medskip
\begin{center}
\begin{tabular}{c|c|c|c|c|c|c|c}
${\texttt{type}}$ &$1^4$ & ${\color{cyan} 1^4}$ & $1\, {\color{cyan} 1}\,2$ & ${\color{cyan} 2^2} $& $1\,3$ & ${\color{cyan}1\,3 }$ & ${\color{cyan} 4}$ \\
\hline   \texttt{size} & $1/48$ & $1/48$ & $1/4$ & $1/8$ & $1/6$ & $1/6$ & $1/4$ \\
\end{tabular} 
\end{center}} 
\caption{{\small The ${\rm H}_4$-conjugacy classes of the elements of ${\rm GL}_2(\Z/3)$.} }\label{tablegl2z3}
\end{table}

In this table, and in others that we will give later, we use standard notations for partitions, and print ${\rm p}_+$ in black and ${\rm p}_-$ in cyan. So the sequence of symbols $1^{a_1}\, {\color{cyan} 1^{b_1}}\,2^{a_2}\,{\color{cyan} 2^{b_2}}\, \dots i^{a_i}\, {\color{cyan} i^{b_i}} \dots$ stands for the couple $(p_+,p_-)$ where $p_+$ is the partition of $\sum_i a_i$ in $a_1$ times $1$, $a_2$ times $2$, and so on, and $p_-$  is the partition of  $\sum_i b_i$ in $b_1$ times $1$, $b_2$ times $2$, and so on. The symbol "$i^m$" (resp. "${\color{cyan}i^{m}}$") is omitted for $m=0$, and replaced by "$i\,$" (resp. "${\color{cyan} i}\,$") for $m=1$. \ps 

\begin{remark} \label{comptablegl2z3} {\rm Table \ref{tablegl2z3} is easily deduced from the conjugacy classes of {\small ${\rm GL}_2(\Z/3)$}. To fix ideas, define the embedding $\iota$ in Example \ref{exgl2z3} by choosing $v_1,v_2,v_3$ and $v_4$ to be respectively {\tiny $\left[\begin{array}{c} 0 \\ 1 \end{array} \right]$, $\left[\begin{array}{c} 1 \\ 0 \end{array} \right]$, $\left[\begin{array}{c} 1 \\ 1 \end{array} \right]$} and {\tiny $\left[\begin{array}{c} 1 \\ -1 \end{array} \right]$}. Then the images under $\iota$ of the elements 
$\pm {\rm I}_2$, {\tiny $\left[\begin{array}{cc} 1 & 0 \\ 0 & -1 \end{array}\right], \left[\begin{array}{cc} 0 & -1 \\ 1 & 0 \end{array}\right], \pm \left[\begin{array}{cc} 1 & 1 \\ 0 & 1 \end{array}\right]$} and {\tiny $\pm \left[\begin{array}{cc} 0 & 1 \\ 1 & 1 \end{array}\right]$} of $\GL_2(\Z/3)$ are respectively {\small $\pm \,1$,  \,\,$\varepsilon_1 \,(34)$,\,\, $\varepsilon_2\,(12)\,\varepsilon_4\,(34)$, \,\,$\pm\, \varepsilon_1 \,\varepsilon_4 \,(134)$ and $\pm\, \varepsilon_2 \,(1342)$}.
}
\end{remark}

\subsection{Digression: subgroups of certain wreath products}\label{parwr}${}^{}$

\indent Let $G$ be a group, $X$ a set equipped with a transitive action of $G$, and $A$ an abelian group. The group $G$ acts in a natural way on the abelian group $A^X$ of all functions $X \rightarrow A$, so we can form the semi-direct  product $H : = A^X \rtimes G$. \ps \ps

We denote by $\pi : H \rightarrow G$ the canonical projection, with kernel ${\rm ker}\, \pi = A^X$.  We have a ``diagonal'' map $\delta : A \rightarrow A^X$, defined by $\delta(a)(x)=a$ for all $a$ in $A$ and $x$ in $X$. This map $\delta$ is an embedding of $G$-modules if we view the source $A$ as a trivial $G$-module: the image of $\delta$ is a central subgroup of $H$. Our aim in this paragraph is to study:\ps\ps

--  the set $\mathcal{C}$ of subgroups $C \subset H$ with $\pi(C)=G$ and ${\rm ker}\, \pi \cap C=1$, \ps\ps

-- the set $\mathcal{G}$ of subgroups $\widetilde{G} \subset H$ with $\pi(\widetilde{G})=G$ and ${\rm ker}\, \pi \cap \widetilde{G} = \delta(A)$. \ps\ps

\noindent The group $H$ acts both on $\mathcal{C}$ and  $\mathcal{G}$ by conjugation. We start with two simple observations: \ps


-- \,For any group morphism $\chi : G \rightarrow A$, the set $G^{\,\chi} := \{ \delta(\chi(g)) \cdot g, \, \, g \in G\}$ is a subgroup of $H$ isomorphic to $G$, and $G^{\,\chi}$ is an element of  $\mathcal{C}$.  \ps\ps

-- \,There is a natural map ${\rm c}_2 :  \mathcal{G} \rightarrow {\rm H}^2(G,A)$, sending $\widetilde{G}$ in $\mathcal{G}$ to the equivalence class of the central extension $1 \rightarrow A \overset{\delta}{\rightarrow} \widetilde{G} \overset{\pi}{\rightarrow} G \rightarrow 1$. Two elements of $\mathcal{G}$ which are $H$-conjugate are also $A^X$-conjugate, hence define the same class in ${\rm H}^2(G,A)$. \ps\ps

We fix some $x \in X$ and denote by $G_x \subset G$ the isotropy group of $x$. For each integer $m \geq 0$, we denote by ${\rm r}_{m} : {\rm H}^m(G,A) \longrightarrow {\rm H}^m(G_x,A)$ the usual restriction map on the cohomology groups of the trivial $G$-module $A$. 

\begin{prop} \label{propwr} \begin{itemize}\item[(i)] For $\chi,\chi' \in {\rm Hom}(G,A)$, the subgroups $G^{\,\chi}$ and 
$G^{\,\chi'}$ of $H$ are conjugate if, and only if, $\chi$ and $\chi'$ coincide on $G_x$. \ps
\item[(ii)]  If ${\rm r}_1$ is surjective then any subgroup $C \in \mathcal{C}$ is conjugate to $G^{\, \chi}$ for some $\chi \in {\rm Hom}(G,A)$. \ps 
\item[(iii)] If ${\rm r}_1$ is surjective then the map ${\rm c}_2 : H\backslash \mathcal{G} \rightarrow {\rm H}^2(G,A)$ is injective, and its image is the subgroup ${\rm ker}\, {\rm r}_2$ of extensions which split over $G_x$. 
\end{itemize}
\end{prop}

\begin{pf} We shall use twice the following classical facts. Let $\Gamma$ be a group acting on an abelian group $V$ and denote by $\pi : V \rtimes \Gamma \rightarrow \Gamma$ the natural projection. Let $\mathcal{K}$ be the set subgroups $K \subset V \rtimes \Gamma$ with $\pi(K)=\Gamma$ and $ {\rm ker}\, \pi\, \cap \,V\,= \,1$. Any $K \in \mathcal{K}$ has the form $\{ s(\gamma) \, \gamma, \, \, \gamma \in \Gamma\}$ for a unique $1$-cocyle $s \in {\rm Z}^1(\Gamma,V)$, that we denote ${\rm s}_K$. The map $K \mapsto {\rm s}_K$, $\mathcal{K} \rightarrow {\rm Z}^1(\Gamma,V)$, is bijective; two elements $K, K'$ in $\mathcal{K}$ are conjugate by an element of $V$ if, and only if, ${\rm s}_K$ and ${\rm s}_{K'}$ have the same class in ${\rm H}^1(\Gamma,V)$. Last but not least, note that $K, K'$ in $\mathcal{K}$ are conjugate by an element of $V$ if, and only if, they are conjugate in $V \rtimes \Gamma$: if we have $K' = g K g^{-1}$ with $g \in V \rtimes \Gamma$, we may write $g = vk$ with $v \in V$ and $k \in K$, and we have $K' = v K v^{-1}$.
\ps\ps

We apply this first to $\Gamma=G$ and $V=A^X$. The map $ {\rm Hom}(G,A)={\rm Z}^1(G,A) \rightarrow {\rm Z}^1(G,A^X)$ defined by $\delta$ sends $\chi$ to the $1$-cocycle defining $G^{\,\chi}$. The choice of $x \in X$ identifies the $G$-module $A^X$ with the co-induced module of the trivial $G_x$-module $A$ to $G$. By Shapiro's lemma, we obtain for each integer $m\geq 0$ a natural isomorphism ${\rm sh} : {\rm H}^m(G,A^X) \lisomo {\rm H}^m(G_x,A)$. Concretely, if $f : G^m \rightarrow A^X$ is an $m$-cocycle, then ${\rm sh}(f)$ is the class of the $m$-cocycle $f' : G_x^m  \rightarrow A$ defined by $f'(g_1,\dots,g_m)=f(g_1,\dots,g_m)(x)$. It follows that the composition of the maps $$ {\rm H}^m(G,A) \overset{{\rm H}^m(\delta)}{\longrightarrow} {\rm H}^m(G,A^X) \overset{{\rm sh}}{\lisomo} {\rm H}^m(G_x,A)$$
coincides with the map ${\rm r}_m$. For $m=1$, this proves assertions (i) and (ii). \ps\ps

Let us prove assertion (iii). Let $Q$ be the cokernel of $\delta$. By applying the first paragraph above to $\Gamma=G$ and $V=Q$, we obtain a natural bijection ${\rm c}_1 : H \backslash \mathcal{G} \lisomo {\rm H}^1(G,Q)$. The long exact sequence of cohomology groups associated to $0 \rightarrow A \overset{\delta}{\rightarrow} A^X \rightarrow Q \rightarrow 0$ contains a piece of the form $$ {\rm H}^1(G,A) \overset{{\rm H}^1(\delta)}{\longrightarrow}  {\rm H}^1(G,A^X) \longrightarrow {\rm H}^1(G,Q) \overset{\eta}{\longrightarrow} {\rm H}^2(G,A)  \overset{{\rm H}^2(\delta)}{\longrightarrow}  {\rm H}^2(G,A^X).$$
By the second paragraph, the kernel of $\eta$ is isomorphic to the cokernel of ${\rm r}_1$, and the image of $\eta$ is the kernel of ${\rm r}_2$. As it is straightforward to check from the definition of ${\rm c}_2$ that we have $\eta \circ c_1 = {\rm c}_2$, this concludes the proof of assertion (iii). 
\end{pf}


\subsection{Applications to ${\rm H}_n$} \label{parhn} ${}^{}$ 
The group ${\rm H}_n$ is of course the special case of the construction of \S \ref{parwr} with $G={\rm S}_n$, $X=\{1,\dots,n\}$ and $A=\{\pm 1\}$ (multiplicative group). The signature $\epsilon$ gives rise to the subgroup ${\rm S}_n^\epsilon$ of ${\rm H}_n$ whose elements have the form $\epsilon(\sigma) \sigma$, $\sigma \in {\rm S}_n$. For any transposition $\tau$ in ${\rm S}_n$ we have ${\rm n}_{-}(\epsilon(\tau) \tau) = n-2$, whereas ${\rm n}_{-}(\sigma) =0$ for all $\sigma$ in ${\rm S}_n$: this shows that ${\rm S}_n^\epsilon$ is not conjugate to ${\rm S}_n$ in ${\rm H}_n$ for $n>2$ (a fact which also follows from assertion (i) below). 

\begin{prop}\label{sectionSn}\label{2Sn} \begin{itemize} \item[(i)] Let $G$ be a subgroup of ${\rm H}_n$ of order $n!$ with $\pi(G)={\rm S}_n$. Then $G$ is either conjugate to ${\rm S}_n$ or to ${\rm S}_n^\epsilon$. Moreover, ${\rm S}_n$ and ${\rm S}_n^\epsilon$ are conjugate in ${\rm H}_n$ if, and only if, we have $n \leq 2$. \ps
\item[(ii)] Let $G$ be a subgroup of ${\rm H}_n$ of order $2 n!$ with $\pi(G)={\rm S}_n$. Then $-1$ is in $G$ and exactly one of the following properties holds: \begin{itemize}
\item[(a)] $G$ is conjugate to $\{ \pm 1\} \cdot {\rm S}_n$, \ps
\item[(b)] $n=2$ and $G \simeq \Z/4$, \ps 
\item[(c)] $n=4$ and $G$ is conjugate to the group ${\rm GL}_2(\Z/3)$ embedded in ${\rm H}_4$ as in Example \ref{exgl2z3}.\end{itemize}\end{itemize}
\end{prop}

\begin{pf} Note first that in case (ii), $\{ \pm 1\}^n \cap G$ is a normal subgroup of order $2$ of $G$, hence it is central and generated by $-1$ by the assumption $\pi(G)={\rm S}_n$. \par
The stabilizer of $n$ in ${\rm S}_n$ is naturally identified with ${\rm S}_{n-1}$, with the convention ${\rm S}_{0}=1$.
The signature $\epsilon$ is a generator of ${\rm H}^1({\rm S}_n,\{\pm 1\})$, so the restriction map ${\rm H}^1({\rm S}_n,\{\pm 1\}) \rightarrow {\rm H}^1({\rm S}_{n-1},\{\pm 1\})$ is clearly surjective, and bijective for $n \neq 2$. Moreover, we know form Schur that the restriction map ${\rm r}_2: {\rm H}^2({\rm S}_n,\{\pm 1\}) \rightarrow {\rm H}^2({\rm S}_{n-1},\{\pm 1\})$ is surjective as well for all $n>1$, and that the dimension of the $\Z/2$-vector space ${\rm H}^2({\rm S}_n,\{\pm 1\})$ is $2$ for $n\geq 4$, $1$ for $n=3$ and $2$, and $0$ for $n=1$ \cite{schur}. The kernel of ${\rm r}_2$ is thus $0$ for $n \neq 2,4$, isomorphic to $\Z/2$ otherwise. We conclude by Proposition \ref{propwr} and Example \ref{exgl2z3}. \end{pf}

\begin{remark} \label{rem:A4} {\rm The natural map ${\rm H}^i({\rm Alt}_4,\Z/2) \rightarrow {\rm H}^i({\rm Alt}_3,\Z/2)$ is $0 \rightarrow 0$ for $i=1$ and $\Z/2 \rightarrow 0$ for $i=2$. By Proposition \ref{propwr} (iii), there is thus a unique conjugacy class of nonsplit central extensions of ${\rm Alt}_4$ by $\{\pm 1\}$ in ${\rm H}_4$ (or in $\{\pm 1\}^4 \rtimes {\rm Alt}_4$). As ${\rm Alt}_4$ does not embed in ${\rm GL}_2(\Z/3)$, one such extension is the inverse image of ${\rm Alt}_4$ in the extension described in Example \ref{exgl2z3}.}
\end{remark}

We now give another example. As is well-known, the group ${\rm S}_5$ has a unique isomorphism class of transitive actions on the set $\{1,\dots,6\}$, obtained from the conjugation action on its $6$ subgroups of order $5$. We fix such an action and consider the associated semi-direct product $\{\pm 1\}^6 \rtimes {\rm S}_5$, as in \S \ref{parwr}. We have a defined {\it loc. cit.} a set $\mathcal{G}$ of subgroups of $\{\pm 1\}^6 \rtimes {\rm S}_5$ which are central extensions of ${\rm S}_5$ by $\{\pm 1\}$.

\begin{prop}\label{corS5} The set $\mathcal{G}$ is the disjoint union of two conjugacy classes: the one of the split extension $\{\pm 1\}\cdot  {\rm S}_5$, and another one consisting of nonsplit extensions which are split over the alternating subgroup ${\rm Alt}_5$ of ${\rm S}_5$. 
\end{prop}

\begin{pf} Let $N \subset {\rm S}_5$ be the normalizer of the subgroup $S=\langle (12345) \rangle$. Then $N$ is the semi-direct product of $\langle (2354) \rangle \simeq \Z/4$ by $S \simeq \Z/5$, so we have ${\rm H}^i(N,\Z/2) \simeq \Z/2$ for each $i\geq 0$ and the restriction map ${\rm H}^1({\rm S}_5,\Z/2) \longrightarrow {\rm H}^1(N,\Z/2)$ is an isomorphism. We observe from the presentation given by Schur of the two Schur-covers of ${\rm S}_5$ that they are non split over the subgroups of ${\rm S}_5$ containing a double transposition, such as $N$ or ${\rm Alt}_5$. This implies that the kernel of the restriction map ${\rm H}^2({\rm S}_5,\Z/2) \longrightarrow {\rm H}^2(N,\Z/2)$ is generated by the remaining nonzero class in ${\rm H}^2({\rm S}_5,\Z/2)$, namely the one which splits over ${\rm Alt}_5$ (recall ${\rm H}^2({\rm Alt}_5,\Z/2) \simeq \Z/2$), and we conclude by Proposition \ref{propwr}. \end{pf}

A homomorphism ${\rm S}_5 \rightarrow {\rm S}_6$ as above can alternatively be constructed from the natural action of ${\rm PGL}_2(\Z/5) \simeq {\rm S}_5$ on the projective line ${\rm P}^1(\Z/5)$. The action of ${\rm GL}_2(\Z/5)$ on the $12$-elements set $((\Z/5)^2-\{0\})/\{\pm 1\}$ permutes the $6$ disjoint pairs of the form $\{v,2v\}$, which defines a natural conjugacy class of embeddings 
\begin{equation} \label{iotaS5} \iota : {\rm GL}_2(\Z/5)/\{\pm {\rm I}_2\} \longrightarrow \{\pm1 \}^6 \rtimes {\rm S}_5.\end{equation} The group $\iota({\rm GL}_2(\Z/5)/\{\pm {\rm I}_2\})$ belongs to the second class of Proposition \ref{corS5} (recall ${\rm PSL}_2(\Z/5) \simeq {\rm Alt}_5$). The map $\iota$ is explicit enough to allow the computation of the conjugacy classes of the elements of $\iota({\rm GL}_2(\Z/5)/\{\pm {\rm I}_2\})$ viewed as a subgroup of ${\rm H}_6 \supset  \{\pm1 \}^6 \rtimes {\rm S}_5$: they are gathered in Table \ref{tab:2S5}.

\begin{table}[htp]
{\scriptsize \renewcommand{\arraystretch}{1.8} \medskip
\begin{center}
\begin{tabular}{c|c|c|c|c|c|c|c|c|c|c|c}
${\texttt{type}}$ & $1^6$ & {\color{cyan}$1^6$} & $1^2 2^2$ & {\color{cyan}$1^2$}$2^2$ & {\color{cyan}$2^3$} & $3^2$ & {\color{cyan}$3^2$} & $1\,${\color{cyan}$1\,$}$4\,$ & $1\,5$ & {\color{cyan}$1\,5$} & {\color{cyan}$6\,$}\\
\hline   \texttt{size} & $1/240$ & $1/240$ & $1/16$ & $1/16$ & $1/12$ & $1/12$ & $1/12$ & $1/4$ & $1/10$ & $1/10$ & $1/6$ 
\end{tabular} 
\end{center}} 
\caption{{\small The ${\rm H}_6$-conjugacy classes of the elements of ${\rm GL}_2(\Z/5)/\{\pm{\rm I}_2\}$.} }\label{tab:2S5}
\end{table}

	
\section{Characteristic masses of root lattices}\label{par:rootlat}

\subsection{Root systems and root lattices} ${}^{}$
Let $V$ be an Euclidean space. By a ${\it root}$ of $V$ we mean an element $\alpha \in V$ with $\alpha \cdot \alpha = 2$; we denote by ${\rm R}(V)$ the set of roots of $V$ (a sphere).  For each $\alpha \in {\rm R}(V)$, the orthogonal reflection about $\alpha$ is an element ${\rm s}_\alpha$ of ${\rm O}(V)$, given by the formula ${\rm s}_\alpha(x)\,=\,x \,-\, (\alpha \cdot x)\, \alpha$. \ps\ps

An {\it ${\bf ADE}$ root system} in $V$ is a finite set $R \subset {\rm R}(V)$ generating $V$ as a real vector space, and such that for all $\alpha,\beta \in R$ we have $\alpha \cdot \beta \in \Z$ and ${\rm s}_\alpha(\beta) \in R$. In particular, $R$ is a root system in the sense of Bourbaki \cite{bourbaki},  and each irreducible component of $R$ is of type ${\bf A}_n$ with $n\geq 1$, ${\bf D}_n$ with $n\geq 4$, or ${\bf E}_n$ with $n=6,7,8$. The {\it root lattice} of $R$ is the lattice ${\rm Q}(R) \overset{{\rm def}}{=} \sum_{\alpha \in R} \Z\, \alpha \subset V$ generated by $R$. This is an even lattice, and we have the important equality
\begin{equation} \label{rootrootlattice} {\rm R}(V) \cap {\rm Q}(R) = R.\end{equation}
\ps

If $L \subset V$ is {\it any integral lattice}, we denote by ${\rm R}(L)=L \cap {\rm R}(V)$ the set of roots of $L$. It follows at once from the definitions that ${\rm R}(L)$ is an ${\bf ADE}$ root system in the Euclidean subspace $U$ of $V$ generated by ${\rm R}(L)$. We say that $L$ is a {\it root lattice} if ${\rm R}(L)$ generates $L$ as an abelian group, {\it i.e.} if we have $L={\rm Q}( {\rm R}(L))$ (hence $U=V$). By definitions and \eqref{rootrootlattice}, the map $R \mapsto {\rm Q}(R)$ is a bijection between the set of ${\bf ADE}$ root systems of $V$ and the set of root lattices of $V$, whose inverse is $L \mapsto {\rm R}(L)$. \ps\ps

 We shall always use a bold font to denote an isomorphism class of root systems, and reserve the normal font for a root lattice with the corresponding root system. For instance, if for $n\geq 2$ we set ${\rm D}_n=\{(x_i) \in \Z^n, \sum_i x_i \equiv 0 \bmod 2\}$ as in the introduction, then ${\rm R}({\rm D}_n)$ is a root system of type ${\bf D}_n$ in the standard Euclidean space $\R^n$. We have also defined {\it loc. cit.} the lattice ${\rm E}_n$ for $n \equiv 0 \bmod 8$. It is easy to check ${\rm R}({\rm E}_n)={\rm R}({\rm D}_n)$ for $n>8$ and that ${\rm R}({\rm E}_8)$ is of type ${\bf E}_8$. We choose in an arbitrary way root lattices ${\rm A}_n$ for $n\geq 1$, as well as ${\rm E}_6$ and ${\rm E}_7$, whose root systems are of type ${\bf A}_n$, ${\bf E}_6$ and ${\bf E}_7$. \ps\ps

Let $L$ be an integral lattice in $V$ and set $R={\rm R}(L)$. The ${\rm s}_\alpha$ with $\alpha$ in $R$ generate a subgroup of ${\rm O}(L)$ called the {\it Weyl group} of $L$, and denoted ${\rm W}(L)$. This is a normal subgroup of ${\rm O}(L)$, and we denote by ${\rm G}(L) = {\rm O}(L)/{\rm W}(L)$ the quotient group.  Assume first that $L$ is the root lattice ${\rm Q}(R)$; in this case we also set ${\rm W}(R) := {\rm W}(L)$, ${\rm O}(R) := {\rm O}(L)$ (this latter group is also denoted ${\rm A}(R)$ by Bourbaki) and ${\rm G}(R)={\rm G}(L)$. 
As is well-known, ${\rm G}(R)$ is isomorphic to the automorphism group of the Dynkin diagram of $R$, and we have
\begin{equation} \label{eq:outdyn} {\rm G}(R) \simeq \left\{ \begin{array}{cl} 1 & \text{for\,\,} R \simeq {\bf A}_1, {\bf E}_7, {\bf E}_8, \\ 
{\rm S}_3 & \text{for\,\,} R \simeq {\bf D}_4, \\
\Z/2 & \text{otherwise.} 
\end{array} 
\right.
\end{equation}
Moreover, ${\rm W}(R)$ permutes the {\it positive root systems}\footnote{Recall that a positive root system in $R$ is a subset of the form $\{ \alpha \in R,\, \varphi(\alpha) >0\}$ where $\varphi : V \rightarrow \R$ is a linear form with  $0 \notin \varphi(R)$.} $R^+$ of $R$, or equivalently the {\it Weyl vectors}\footnote{A Weyl vector of $R$ is a vector of the form $\rho = \frac{1}{2} \sum_{\alpha  \in R^+} \alpha$ for $R^+$ a positive root system of $R$. In particular we have $2\rho \in {\rm Q}(R)$.}of $R$, in a simply transitive way. Let us now go back to the case of an arbitrary $L$. The set $R$ is a root system in the Euclidean space $U$ generated by $R$, and the restriction $\sigma \mapsto \sigma_{|U}$ induces a morphism ${\rm O}(L) \rightarrow {\rm O}(R)$ and an isomorphism  ${\rm W}(L) \isomo {\rm W}(R)$. It follows that ${\rm O}(L)$ permutes the Weyl vectors of $R$, and that ${\rm W}(L)$ permutes them simply transitively. So for any Weyl vector $\rho$ of $R$, the stabilizer ${\rm O}(L)_\rho$ of $\rho$ in ${\rm O}(L)$ is naturally isomorphic to ${\rm G}(L)$ and we have
\begin{equation}\label{splitseqOL} {\rm W}(L) \cap {\rm O}(L)_\rho = 1, \hspace{.3cm} {\rm O}(L)\, =\, {\rm W}(L) \cdot {\rm O}(L)_\rho \hspace{.3cm} \text{and}  \hspace{.3cm} {\rm W}(L) \simeq {\rm W}(R),\end{equation} 
so that ${\rm O}(L)$ is the semi-direct product of ${\rm O}(L)_\rho$ by ${\rm W}(L)$. 

\subsection{Characteristic masses of irreducible root lattices}${}^{}$\label{par:irrmasses}

In this paragraph, we let $R \subset V$ be an ${\bf ADE}$ root system and $L={\rm Q}(R)$. Set $n=\dim V$. Our aim is to determine the characteristic masses of ${\rm O}(L)$ and, more generally, the map ${\rm m}_S : {\rm Car}_n \rightarrow \Q_{\geq 0}$ where $S$ is any subset of the form $\sigma {\rm W}(L)$ with $\sigma \in {\rm O}(L)$ (see \S \ref{par:introdiminv} (ii) for the definition of ${\rm m}_S$). We assume first $R$ irreducible, and argue case by case. \ps\ps

\ps \ps
{\sc (A)\, Case $R \simeq {\bf A}_n$ with $n\geq 1$.} 
\ps\ps

We may assume that $V$ is the hyperplane of sum $0$ vectors in $\R^{n+1}$ and $R=\{\pm (\epsilon_i - \epsilon_j), \, \, 1 \leq i < j \leq n+1\}$, where $\epsilon_1,\dots,\epsilon_{n+1}$ denotes the canonical basis of $\R^{n+1}$, and $L={\rm A}_n$. The group ${\rm W}({\rm A}_n)$ may be identified with the symmetric group ${\rm S}_{n+1}$, acting on $V$ by permuting coordinates.\ps\ps
\newcommand{\um}{\underline{m}}
Let $\mathcal{S}$ denote the set of integer sequences $\um=(m_i)_{i \geq 1}$ with $m_i \geq 0$ for each $i$, and $m_i=0$ for $i$ big enough. Let $\mathcal{A}_{n} \subset \mathcal{S}$ denote the subset of $\underline{m}$ such that  $\sum_i \, i \, m_i \,=\, n+1$. For any $\underline{m}$ in $\mathcal{A}_{n}$, the elements of ${\rm S}_{n+1}$ whose cycle decomposition contains $m_i$ cycles of length $i$ for each $i$ form a single conjugacy class ${\rm C}_{\underline{m}} \subset {\rm S}_{n+1}$. We have furthermore 
$|{\rm C}_{\underline{m}}|=(n+1)!/{\rm n}_{\underline{m}}$ with 
$${\rm n}_{\underline{m}} \,=\, \prod_i \, m_i !\,\, i^{m_i}.$$ 
The characteristic polynomial of ${\rm C}_{\underline{m}}$ acting on $V$ is 
$${\rm P}_{\underline{m}}  \,=\, (t-1)^{-1}\, \prod_i \,(t^i -1)^{m_i},$$
since $\R^{n+1}/V$ is the trivial representation of ${\rm S}_{n+1}$. The following trivial lemma even shows that we have ${\rm P}_{\underline{m}} \neq {\rm P}_{\underline{m'}}$ for $\underline{m} \neq \underline{m'}$. 

\begin{lemma} \label{uniquenesscyc}The polynomials $t^l-1$, with $l\geq 1$, are $\Z$-linearly independent in the multiplicative group of the field $\Q(t)$. \end{lemma}

As a consequence, we obtain the:

\begin{cor} For $n\geq 1$, we have ${\rm m}_{{\rm W}({\rm A}_n)}({\rm P}_{\underline{m}}) \,= \,1/{\rm n}_{\underline{m}}$ for $\underline{m}$ in $\mathcal{A}_{n}$, and ${\rm m}_{{\rm W}({\rm A}_n)}(P)=0$ for all other $P$ in ${\rm Car}_n$.
\end{cor}

As is easily seen, the element $-1= - {\rm id}_V$ is in ${\rm W}({\rm A}_n)$ if and only if $n=1$, and we have ${\rm O}({\rm A}_n) = {\rm W}({\rm A}_n) \cup - {\rm W}({\rm A}_n)$ (this fits of course Formula \eqref{eq:outdyn}). The map ${\rm m}_S$ for the coset $S = - {\rm W}({\rm A}_n)$ is deduced from ${\rm m}_{{\rm W}({\rm A}_n)}$ by the following trivial lemma: 

\begin{lemma}\label{mmoinsS}  Let $S$ be a finite subset of ${\rm O}(V)$ with $\dim V=n$. Then for all $P \in \R[t]$ we have ${\rm m}_{-S}(P) = {\rm m}_{S}(Q)$ with $Q(t) = (-1)^nP(-t)$.
\end{lemma}

\ps \ps
{\sc (D)\, Case $R \simeq {\bf D}_n$ with $n\geq 3$.}
\ps\ps
 We may assume $V = \R^n$, $R = \{ \pm \epsilon_i  \pm \epsilon_j, \, \, 1 \leq i < j \leq n\}$ where $\epsilon_1,\dots,\epsilon_n$ denote again the canonical basis of $V$, and $L={\rm D}_n$. The lattice ${\rm D}_n$ is the largest even sublattice of the standard lattice 
	$${\rm I}_n = \Z^n = \oplus_{i=1}^n \Z \epsilon_i,$$
and thus ${\rm O}({\rm I}_n)$ is a subgroup of ${\rm O}({\rm D}_n)$. This group ${\rm O}({\rm I}_n)$ is nothing else than the hyperoctahedral group ${\rm H}_n$ already introduced in \S \ref{par:defhn}: we have 
$${\rm O}({\rm I}_n) = {\rm H}_n = \{\pm 1\}^n \rtimes {\rm S}_n$$
where ${\rm S}_n$ (resp. $\{ \pm 1\}^n$) act on $\R^n$ by permuting coordinates (resp. sign changes). As is well-known, ${\rm W}({\rm D}_n)$ is the index $2$ subgroup $\ker\,{\rm s}$ of ${\rm O}({\rm I}_n)$ (recall ${\rm s}$ is defined by Formula \eqref{defshn}). By \eqref{eq:outdyn} we also have 
\begin{equation} {\rm O}({\rm I}_n)={\rm O}({\rm D}_n) \, \, \, \, \text{for}\,\, n \neq 4\,\,\, \, \, \, \text{and}\, \, \,{\rm G}({\rm D}_4) \simeq {\rm S}_3 \,\,\,\text{(triality).}\end{equation}
\newcommand{\ump}{\underline{m}^+}
\newcommand{\umm}{\underline{m}^-}
The conjugacy classes of ${\rm H}_n$ have been recalled in \S \ref{par:conjhn}. Let $\mathcal{D}_n \subset \mathcal{S} \times \mathcal{S}$ be the subset of $(\ump,\umm)$ with $\sum_i \, i\, (m^+_i + m^-_i)\, =n$. For any $(\ump,\umm)$ in $\mathcal{D}_n$ the elements of ${\rm H}_n$ whose cycle decomposition contains $m_i^+$ (resp. $m_i^-$) cycles of length $i$ with sign $+1$ (resp. $-1$) for each $i$  form a single conjugacy class ${\rm C}_{\ump,\umm} \subset {\rm H}_{n}$. We easily check  
$|{\rm C}_{\ump,\umm}|=2^n n! / {\rm n}_{\ump,\umm}$ with 
$${\rm n}_{\ump,\umm} \,=\, \prod_i \,\, m^+_i!\,\, \,m^-_i!\, \,\,(2i)^{m^+_i+m^-_i},$$ 
and ${\rm s}({\rm C}_{\ump,\umm}) = (-1)^{|\umm|}$ where we have set $|\um|=\sum_i m_i$ for $\um \in \mathcal{S}$. 
\newcommand{\unp}{\underline{n}^+}
\newcommand{\unm}{\underline{n}^-}
The characteristic polynomial of ${\rm C}_{\ump,\umm}$ acting on $V$ is 
$${\rm P}_{\ump,\umm}  \,=\, \prod_i \,(t^i -1)^{m^+_i}\,(t^i+1)^{m^-_i}= \prod_i \,(t^i -1)^{m^+_i-m^-_i+m^-_{i/2}},$$
\noindent where we have set $m^-_{i/2}=0$ for $i$ odd, and used for $i\geq 1$ the relation $(t^i-1)(t^i+1)=(t^{2i}-1)$. In contrast with the ${\bf A}_n$ case, we may thus have ${\rm P}_{\ump,\umm}={\rm P}_{\unp,\unm}$ for distinct $(\ump,\umm)$ and $(\unp,\unm)$ in $\mathcal{D}_n$. This leads us to introduce the subset 
$$\mathcal{D}_n' = \{ (\ump,\umm) \in \mathcal{D}_n \, \,|\,\, m_i^+m_i^-=0\, \, \, \text{for all}\, i\geq 1\}.$$
Lemma  \ref{uniquenesscyc} shows that we have ${\rm P}_{\ump,\umm} \neq {\rm P}_{\unp,\unm}$ for $(\ump,\umm) \neq (\unp,\unm)$ in $\mathcal{D}_n'$. We reduce to $\mathcal{D}_n'
$ as follows. Consider the following map $\phi : \mathcal{D}_n \rightarrow \mathcal{D}_n$:

\ps\ps

(i) if $(\ump,\umm) \in \mathcal{D}'_n$ set $\phi(\ump,\umm)=(\ump,\umm)$,\ps\ps 

(ii) otherwise there is a smallest $j\geq 1$ with $m_j^+m_j^- \neq 0$ and we set $\phi(\ump,\umm)=(\unp,\unm)$ with  $(n_i^+,n_i^-)=(m_i^+,m_i^-)$ for $i\neq j$ or $i \neq 2j$, and with $(n_j^+,n_j^-)=(m_j^+-1,n_j^+-1)$ and $(n_{2j}^+,n_{2j}^-)=(m_{2j}^++1,m_{2j}^-)$.\ps\ps

\noindent It is clear that we have ${\rm P}_{\phi(\ump,\umm)}={\rm P}_{\ump,\umm}$ for all $(\ump,\umm)$ in $\mathcal{D}_n$, and that for each $m=(\ump,\umm) \in \mathcal{D}_n$ the sequence $m,\phi(m),\phi^2(m),\dots$ is eventually constant and equal to some element of $\mathcal{D}_n'$, that we denote by $\psi(m)$. 

\begin{cor}\label{car:massesDn} Let $\sigma \in {\rm O}({\rm I}_n)$. For all $(\ump,\umm)$ in $\mathcal{D}_n'$ we have 
$${\rm m}_{\sigma {\rm W}({\rm D}_n)}({\rm P}_{\ump,\umm}) \,=\, \sum\,\,\frac{1}{{\rm n}_{\unp,\unm}}$$
the sum being over all the $(\unp,\unm)$ in $\mathcal{D}_n$ with $\psi(\unp,\unm)=(\ump,\umm)$ and 
$(-1)^{|\umm|} = {\rm s}(\sigma)$. We have ${\rm m}_{\sigma {\rm W}({\rm D}_n)}(P)=0$ for all other $P$ in ${\rm Car}_n$.
\end{cor}

We have ${\rm G}({\rm D}_4) \simeq {\rm S}_3$ so it remains to determine ${\rm m}_{\sigma {\rm W}({\rm D}_4)}$ for the $6$ possible classes $\sigma {\rm W}({\rm D}_4)$. A first general reduction is the following lemma:

\begin{lemma}\label{lem:conjcoset} Let $L$ be an integral lattice, as well as elements  $\sigma_1,\sigma_2$ in ${\rm O}(L)$ whose images in ${\rm G}(L)$ are conjugate. Then we have ${\rm m}_{\sigma_1{\rm W}(L)}={\rm m}_{\sigma_2{\rm W}(L)}$. 
\end{lemma}

\begin{pf} Write $\sigma_2 = \gamma \sigma_1 \gamma^{-1} w_0$ with $\gamma$ in ${\rm O}(L)$ and $w_0$ in ${\rm W}(L)$.  For $w \in {\rm W}(L)$ we have 
$\det (t - \sigma_2 w) = \det (t - \sigma_1 \gamma^{-1} w_0 w \gamma)$.  We conclude as $w \mapsto  \gamma^{-1} w_0 w \gamma$ is a bijection of the normal subgroup ${\rm W}(L)$ of ${\rm O}(L)$. \end{pf}

In particular, ${\rm m}_{\sigma {\rm W}({\rm D}_4)}$ is already given by Lemma \ref{car:massesDn} whenever the image of $\sigma$ in ${\rm G}({\rm D}_4) \simeq {\rm S}_3$ has order $1$ or $2$, and does not depend on $\sigma$ if this image has order $3$. There are many ways to determine ${\rm m}_{\sigma {\rm W}({\rm D}_4)}$ in this latter case. One way is to consider first the set 
$$R'=\{ v \in {\rm D}_4, v \cdot v = 4\} = \{ \pm 2 \epsilon_i\, \, | \, \, i =1,\dots,4\} \cup \{ \sum_{i=1}^4 \pm \epsilon_i\}$$
and observe that we have $\alpha \cdot x \in 2 \Z$ for all $\alpha \in R'$ and $x \in {\rm D}_4$. In particular, $\frac{1}{\sqrt{2}} R'$ is a root system (of type ${\bf D}_4$) in $\R^4$ and we have ${\rm W}(\frac{1}{\sqrt{2}} R') \subset {\rm O}({\rm D}_4)$. The two roots $\alpha = \sqrt{2} \epsilon_1$ and $\beta = \frac{1}{\sqrt{2}} (\epsilon_1 + \epsilon_2 + \epsilon_3 + \epsilon_4)$ are in $\frac{1}{\sqrt{2}} R'$ with $\alpha \cdot \beta =1$, and the order $3$ element 
$$\sigma_0 :={\rm s}_\beta \circ {\rm s}_\alpha\, =\,  {\tiny \,\frac{1}{2}\, \left(\begin{array}{cccc} -1 & -1 & -1 & -1 \\ 1 & 1 & -1 & -1 \\ 1 & -1 & 1 & -1 \\ 1 & -1 & -1 & 1 \end{array}\right)}$$
\noindent of ${\rm O}({\rm D}_4)$ does not belong to ${\rm W}({\rm D}_4)$. It is trivial to enumerate with a computer the $2^3 \,4! = 384$ elements of $\sigma_0 {\rm W}({\rm D}_4)$ and to list their characteristic polynomials.  We obtain:

\begin{cor} Let $\sigma$ be an element of order $3$ in ${\rm G}({\rm D}_4)$. The ${\rm m}_{\sigma {\rm W}({\rm D}_4)}(P)$ with $P \in {\rm Car}_4$ are given by Table\footnote{The notation for polynomials used in Table \ref{TableCWD4} is explained on p.~ \pageref{not:tables}.} \ref{TableCWD4}.  
\end{cor} 

In fact, the reasoning above can be pushed a little further: it turns out that $R''=R \coprod R'$ is a root system of type ${\bf F}_4$ in $\R^4$ (not ${\bf ADE}$ of course) and that we have ${\rm W}(R'') = \langle {\rm W}(R), {\rm W}(\frac{1}{\sqrt{2}} R') \rangle = {\rm O}({\rm D}_4)$. But the conjugacy classes of ${\rm W}(S)$, with $S$ any irreducible root system of exceptional type, have been listed and studied in a conceptual way by Carter in \cite{carter}, including their characteristic polynomials (see p. 22 \& 23 {\it loc. cit.}).  The map ${\rm m}_{{\rm O}({\rm D}_4)}$ may be deduced in particular from Table 8 of  \cite{carter}. The map ${\rm m}_{\sigma_0 {\rm W}({\rm D}_4)}$ follows then from the equality
${\rm m}_{{\rm O}({\rm D}_4)} \,=\,- \frac{1}{3} \,{\rm m}_{{\rm W}({\rm D}_4)} \,+\, \frac{1}{2} \,{\rm m}_{{\rm O}({\rm I}_4)}\, +\, \frac{1}{3}\, {\rm m}_{\sigma_0 {\rm W}({\rm D}_4)}$.\ps\ps

\begin{remark}\label{rem:ordertauwr}
{\rm
Assume $R$ is an irreducible root system. It follows from \eqref{eq:outdyn} that two elements of ${\rm G}(R)$ are conjugate if and only if they have the same order, which is always $1, 2$ or $3$. In particular, Lemma \ref{lem:conjcoset} shows that for $\sigma$ in ${\rm O}(R)$ the map ${\rm m}_{\sigma {\rm W}(R)}$ only depends on the order of $\sigma$  in ${\rm G}(R)$.
}
\end{remark}

\ps \ps
{\sc (E)\, Cases $R \simeq {\bf E}_n$ with $n=6, 7$ and $8$.}
\ps\ps

The aforementioned results of Carter also allow to deduce ${\rm m}_{\rm W({\rm E}_n)}$ for $n=6,7$ and $8$ (using Tables 9, 10 and 11 {\it loc. cit}). Alternatively, and as a useful check, these masses can also be computed directly using a variant of the {\sff Algorithm A} explained in \S \ref{par:introalgo}. Indeed, choosing a positive system $R^+ \subset R$, we may view ${\rm W}(R)$ as the subgroup of ${\rm O}(V)$ generated by the $n$ reflections ${\rm s}_\alpha$, with $\alpha$ a simple root in $R^+$. As ${\rm W}(R)$ acts faithfully and transitively on $R$, it is also the subgroup of the permutation group of $R$ generated by these $n$ permutations ${\rm s}_\alpha$, with $|R|=72$ (case $n=6$), $|R|=126$ (case $n=7$) or $|R|=240$ (case $n=8$). Applying \texttt{GAP}'s \texttt{ConjugacyClasses} algorithm to this permutation group, we obtain representatives and cardinalities of the conjugacy classes of ${\rm W}(R)$, and it only remains to compute their characteristic polynomials. All in all, these computations only take a few seconds for the computer.
Both methods lead to the:

\begin{cor} For $n=6,7$ and $8$, the ${\rm m}_{{\rm W}({\rm E}_n)}(P)$ with $P \in {\rm Car}_n$ are given by Tables \ref{TableWE6}, \ref{TableWE7} and \ref{TableWE8}.  
\end{cor}

Note that for $n=7,8$ we have ${\rm O}({\rm E}_n) = {\rm W}({\rm E}_n)$ (no non trivial diagram automorphism). For $n=6$, we have ${\rm O}({\rm E}_6) = {\rm W}({\rm E}_6) \coprod - {\rm W}({\rm E}_6)$, but the map ${\rm m}_{-{\rm W}({\rm E}_6)}$ is deduced from ${\rm m}_{{\rm W}({\rm E}_6)}$ using Lemma \ref{mmoinsS}.\ps\ps

\subsection{The non irreducible case} \label{par:nonirr}${}^{}$

Assume now $R$ is a non necessarily irreducible ${\bf ADE}$ root system in $V$, set $L={\rm Q}(R)$ and fix $\sigma$ in ${\rm O}(R)$. Our aim is to give a formula for ${\rm m}_{\sigma {\rm W}(R)}$. Write $R$ as the disjoint union of its irreducible components $ R = \coprod_{i \in I} R_i$. We have
$$L = \overset{\perp}{\bigoplus_{i \in I}}\,\,{\rm Q}(R_i) \, \,\, \text{and} \,\,\, {\rm W}(R) = \prod_{i \in I} {\rm W}(R_i).$$ The element $\sigma$ induces a permutation of the set $\{ R_i \,|\, \, \, i \in I\}$ of irreducible components of $R$. We write $\sigma \,=\, c_1\,c_2 \cdots c_r$ the cycle decomposition of this permutation. For each $j=1,\dots,r$, we choose an irreducible component $S_j$ of $R$ in the support of $c_j$, denote by $s_j = \dim {\rm Q}(S_j)$ the rank of $S_j$ and by $l_j$ the length of the cycle $c_j$. For each $j$ we have $\sigma^{l_j}(S_j)=S_j$ and we denote by $\tau_j$ the restriction of $\sigma^{l_j}$ to ${\rm Q}(S_j)$;  so $\tau_j$ is an element of ${\rm O}(S_j)$.\ps\ps

\begin{prop} \label{prop:massnonirr}
In the setting above, we have for all $P$ in ${\rm Car}_n$  
$${\rm m}_{\sigma {\rm W}(R)}(P) = \sum_{(P_1,\dots,P_r)} \,\,\, \, \prod_{j=1}^r {\rm m}_{\tau_j {\rm W}(S_j)}(P_j)$$
summing over all $(P_1,\dots,P_r) \in {\rm Car}_{s_1} \times \cdots \times {\rm Car}_{s_r}$ with $\prod_{j=1}^r P_j(t^{l_j})=P(t)$.
\end{prop}

The first ingredient in the proof is the following trivial lemma. 

\begin{lemma}
\label{lem:prodmasses}
For $i=1,2$, let $V_i$ be an Euclidean space and $\Gamma_i \subset {\rm O}(V_i)$ a finite subset. Set $V=V_1 \perp V_2$ and view $\Gamma = \Gamma_1 \times \Gamma_2$ as a subset of ${\rm O}(V)$. For all monic polynomials $P$ in $\R[t]$ of degree $\dim V$ we have $${\rm m}_\Gamma(P) \,=\, \sum_{(P_1,P_2)}  \, {\rm m}_{\Gamma_1}(P_1)\, \, {\rm m}_{\Gamma_2}(P_2),$$
the sum being over the $(P_1,P_2)$, with $P_i \in \R[t]$ monic of degree $\dim V_i$, and with $P_1P_2=P$.
\end{lemma}

\begin{pf} (of Proposition \ref{prop:massnonirr}) Applying Lemma \ref{lem:prodmasses}, we may and do assume $r=1$, {\it i.e.} that $\sigma$ permutes transitively the irreducible components of $R$. In this case, we simply write $(S,s,l,\tau)$ instead of $(S_1,s_1,l_1,\tau_1)$.
We may also assume that we have $I=\{0,\dots,l-1\}$ and $R_i=\sigma^{i}(S)$ for $0 \leq i < l$. In particular, we have 
${\rm W}(R)\, =\, \prod_{0 \leq i < l}\,\,\,\sigma^i \,{\rm W}(S)\, \sigma^{-i}$.  Choose a $\Z$-basis $e=(e_1,\dots,e_s)$ of ${\rm Q}(S)$ and consider the following $\Z$-basis of $L$: $$f=(e_1,\dots,e_s,\sigma(e_1),\dots,\sigma(e_s),\dots,\sigma^{l-1}(e_1),\dots,\sigma^{l-1}(e_s)).$$ For all $w=(\sigma^i w_i \sigma^{-i})_{0 \leq i <l}$ in ${\rm W}(R)$, the matrix of $\sigma w$ in the basis $f$ is
{\small \begin{equation} \label{eq:matform}\left( \begin{array}{ccccc}  &  &  &  & M W_{l-1} \\ W_0 &  &  &  &  \\   & W_1 &  &  &  \\  & & \ddots & &  \\  &  &  & W_{l-2} &  \end{array}\right) \end{equation}}

\noindent where $W_i$ is the matrix of  $w_i \in {\rm W}(S)$ in the basis $e$, and $M$ is the matrix of $\tau=\sigma^l$ in the basis $e$. By Lemma \ref{lem:cayleyham} below, it follows that the multiset of polynomials $\det (t - \sigma w)$ (counted with their multiplicities) when $w$ varies in ${\rm W}(R)$, coincides with that of polynomials $\det(t^s -\tau  w_{l-1} w_{l-2} \cdots w_0)$ when the $l$-tuple $(w_0,w_1,\dots,w_{l-1})$ varies in ${\rm W}(S)^l$. As ${\rm W}(S)$ is a group, this multiset is also $|{\rm W}(S)|^{l-1}$ times the multiset of the $\det (t^s- \tau w)$ when $w$ varies in ${\rm W}(S)$, and we are done.
\end{pf}

\begin{lemma} \label{lem:cayleyham} For any $W_0,W_1,\dots,W_{l-1}$ and $M$ in ${\rm M}_s(\C)$, the characteristic polynomial of the matrix \eqref{eq:matform} of size $sl$ is $\det(t^s- M W_{l-1} W_{l-2} \cdots W_0)$. 
\end{lemma}

\begin{pf} By continuity, we may assume $W_i \in {\rm GL}_s(\C)$ for each $i$. Up to conjugating \eqref{eq:matform} by the diagonal matrix 
$(1,W_0,W_1W_0,\dots,W_{l-2} \cdots W_1 W_{0})$
we may assume $W_0=W_1=\cdots = W_{l-1}= {\rm id}_s$. But in this case, the entries of \eqref{eq:matform} commute and we conclude by \cite{ing}
and the following well-known fact  (applied to $Q=t^s-a$): the characteristic polynomial of the companion matrix of a given monic polynomial $Q$ is the polynomial $Q$ itself. 
\end{pf} 

\section{An algorithm computing characteristic masses} \label{par:genalgo}

\subsection{{\sff Algorithm B}} ${}^{}$ \label{par:quickalgo}

Consider the following algorithm, which takes as input an integral lattice $L$ in the standard Euclidean space $V=\R^n$:
{\sff
\small
\ps\ps
B1. Compute the root system $R = {\rm R}(L)$, a positive root system $R^+ \subset R$ and the associate Weyl vector $\rho = \frac{1}{2}\sum_{\alpha \in R^+} \alpha$.\ps
B2. Determine a set $\mathcal{G}$ of generators of the stabilizer ${\rm O}(L)_\rho$ of $\rho$ in ${\rm O}(L)$.\ps
B3. Compute the set $\mathcal{S} = {\rm S}(L)$ defined on p. \pageref{def:S(L)} and view ${\rm O}(L)_\rho$ as the subgroup of permutations of $\mathcal{S}$ generated by $\mathcal{G}$.\ps
B4. Use permutation groups algorithms to determine the sizes $(m_j)_{j \in J}$ and representatives $(\gamma_j)_{j \in J}$ of the conjugacy classes $(c_j)_{j \in J}$ of ${\rm O}(L)_\rho$. \ps
B5. Compute the set ${\rm Irr}(R)$ of irreducible components of $R$, the isomorphism class of each such component, and a basis of the orthogonal $R^\perp$ of $R$ in $V$.\ps
B6. For each $j$ in $J$, compute: \ps
\hspace{.3cm}-- the characteristic polynomial $P_j$ of $\gamma_j$ on $R^\perp$, \ps
\hspace{.3cm}-- a set of representatives ${\rm Irr}_j \subset {\rm Irr}(R)$ of the orbits for the action of $\gamma_j$ on ${\rm Irr}(R)$, \ps
\hspace{.3cm}-- for each $S \in {\rm Irr}_j$, the size $l_S$ of its $\gamma_j$-orbit and the order $d_S \in \{1,2,3\}$ of the permutation $\gamma_j^{l_S}$ of $S$.  \ps
B7. For each $(S,d_S) \in {\rm Irr}(R) \times \{1,2,3\}$ found in B6, compute ${\rm m}_{\tau {\rm W}(S)}$ using the results of \S \ref{par:irrmasses}, where $\tau$ is any element of order $d_S$ in ${\rm O}(S)/{\rm W}(S)$ (see Remark \ref{rem:ordertauwr}). \ps
B8. Using Proposition \ref{prop:massnonirr} and step B7, deduce for each $j$ in $J$ the map ${\rm m}_{\gamma_j {\rm W}(R)}$. \ps
B9. For each $j$ in $J$, define $M_j : {\rm Car}_n \rightarrow \Q$ by setting $M_j(P)={\rm m}_{\gamma_j {\rm W}(R)}(Q)$ if we have $P=QP_j$, and $M_j(P)=0$ otherwise.\ps
B10. Return $\frac{\sum_{j \in J} m_j \, M_j}{\sum_{j \in J} m_j}$. \ps
}
\ps\ps

We will say more about each step of this algorithm in \S \ref{par:implementalgo}. Recall from \eqref{splitseqOL} that we have a semi-direct product ${\rm O}(L) = {\rm W}(L) \rtimes {\rm O}_\rho(L)$ and that the restriction to the subspace $U= {\rm Q}(R) \otimes \R$ of $V$ induces a morphism ${\rm res}: {\rm O}(L) \rightarrow {\rm O}(R)$, an isomorphism ${\rm W}(L) \isomo {\rm W}(R)$ and a morphism ${\rm O}(L)_\rho \rightarrow {\rm O}(R)_\rho$. Together with \eqref{eq:outdyn},  this explains why the elements $d_S$ introduced in the step {\sff B6} are indeed in $\{1,2,3\}$. Moreover, the more correct notation for $\gamma_j {\rm W}(R)$ in {\sff B8} should be ${\rm res}(\gamma_j)  {\rm W}(R)$. For $j$ in $J$, we have 
$M_j={\rm m}_{\gamma_j {\rm W}(L)}$ as ${\rm W}(L)$ acts trivially on $R^\perp$. Last but not least, Lemma \ref{lem:conjcoset} shows $${\rm m}_{{\rm O}(L)} = \frac{\sum_{j \in J} m_j \,\,{\rm m}_{\gamma_j {\rm W}(L)}}{\sum_{j \in J} m_j}.$$ We have proved the:

\begin{prop} {\rm \sff Algorithm B} returns ${\rm m}_{{\rm O}(L)}$. \end{prop}

\subsection{Precisions and an implementation}\label{par:implementalgo}${}^{}$

We now discuss more precisely the steps of {\sff Algorithm B}, as well as some aspects of our implementation: see \cite{homepage} for the source code and a documentation of the \texttt{PARI/GP} function \texttt{masses\char`_calc} (requiring \texttt{GAP}) that we developped. Its input is a Gram matrix \texttt{G} of the lattice $L$, which is thus viewed as the lattice $\Z^n$ equipped with the inner product defined by \texttt{G}. \ps\ps

{\sff B1}. Apply the Fincke-Pohst algorithm \cite{FP} to \texttt{G} to compute $R \subset \Z^n$. In \texttt{PARI}'s implementation, \texttt{qfminim(G)[3]} returns a set $T \subset \Z^n$ with $T \cup -T = R$ consisting of all the elements of $R$ lying in a certain half-space of $\R^n$: this is a positive system, and we simply choose $R^+=T$. \ps\ps

{\sff B2}.  Let \texttt{b} be the Gram matrix of the $\Z$-valued bilinear form $(x,y) \mapsto \,4 (\rho \cdot x)(\rho \cdot y)$ in the canonical basis of $\Z^n$. Apply the Plesken-Souvignier algorithm~\cite{pleskensouvignier}  to the pair of matrices $(\texttt{G},\texttt{b})$. This is implemented in \texttt{PARI/GP} as $\texttt{qfauto([G,b])[2]}$ (following Souvignier's \texttt{C} code). It returns a set $\mathcal{G}' \subset {\rm GL}_n(\Z)$ of generators of the subgroup of ${\rm O}(L)$ whose elements $g$ satisfy $g \rho = \pm \rho$. For each $g \in \mathcal{G}'$, determine the sign $\epsilon_g$ with $g \rho = \epsilon_g \rho$. Define\footnote{An alternative (cleaner) method to compute the stabilizer in ${\rm O}(L)$ of a given element $x$ of $L$ would be to simply add the condition $v_i \cdot x = b_i \cdot x$ for all $i \leq k$ in the definition of a $k$-partial automorphism in \S 3 of \cite{pleskensouvignier}, as well as a similar constraint in the definition of their fingerprint in \S 4 {\it loc. cit}. 
The main advantage of the trick we use is that we do not have to modify the code of the \texttt{PARI} port of Souvignier's program.} $\mathcal{G}$ as the set of $\epsilon_g g$ with $g \in \mathcal{G}'$. \ps\ps

{\sff B3}. Apply recursively the Fincke-Pohst algorithm to find $\mathcal{S} \subset \Z^n$, as explained on p. \pageref{def:S(L)}.
Choose arbitrarily an ordering $\psi : \mathcal{S} \isomo \{1,\dots,N\}$. For each $g$ in $\mathcal{G}$, compute the permutation $\sigma_g = \psi \circ g \circ \psi^{-1}$ in the symmetric group ${\rm S}_{N}$. 
For later use, also extract a basis $\mathcal{S}_0 \subset \mathcal{S}$ of $\R^n$. \ps \ps

{\sff B4}. Apply \texttt{GAP}'s \texttt{ConjugacyClasses} algorithm to the subgroup $H$ of ${\rm S}_{N}$ generated by the $\sigma_g$ with $g$ in $\mathcal{G}$. It returns a list of representatives $(r_j)_{j \in J}$ of the conjugacy classes of $H$, as well as their cardinalities $(m_j)_{j \in J}$. Each $r_j$ is a permutation of $\{1,\dots,N\}$. Using the subset $\mathcal{S}_0$ introduced in {\sff B3}, compute the matrix $\gamma_j \in {\rm GL}_n(\Z)$ of the element of ${\rm O}(L)_\rho$ corresponding to $r_j$ under the natural isomorphism $H \simeq {\rm O}(L)_\rho$. \ps\ps

{\sff B5}. Compute first the basis $B$ of the root system $R$ associated to $R^+$, using $B=\{ \alpha \in R^+ \, \, |\, \, \alpha \cdot \rho =1\}$. Define a graph with set of vertices $B$, and with an edge between $b,b' \in B$ if and only if we have $b \cdot b' \neq 0$. Determine the connected components $B = \coprod_{i \in I} B_i$ of this graph. For $i$ in $I$ define $R_i^+$ as the subset of elements $\alpha$ in $R^+$ with $\alpha \cdot B_i \neq 0$. We have ${\rm Irr}(R) = \{ R_i \, \, |\, \, i \in I\}$. The isomorphism class of the ${\bf ADE}$ root system $R_i = R_i^+ \cup -R_i^+$ is uniquely determined by its rank $|B_i|$ and its cardinality $2|R_i^+|$. \ps\ps

{\sff B6}. Use $i \mapsto R_i$ to identify $I$ with ${\rm Irr}(R)$. Compute the Weyl vector $\rho_i = \frac{1}{2} \sum_{\alpha \in R_i^+} \alpha$ of $R_i$ for each $i$ in $I$.
Fix $j \in J$. There is a unique permutation $\tau_j$ of $I$ such that $\gamma_j(\rho_i)=\rho_{\tau_j(i)}$ for all $i$ in $I$. Compute $\tau_j$
and determine its cycle decomposition. \ps\ps

Steps {\sff B7-B10} are theoretically straightforward. Nevertheless, the efficient implementation of these steps depends on the way the maps ${\rm m}_S$ are represented: see the documentation in \cite{homepage} for more about the (imperfect) way we proceed in \texttt{gp}. In the end,  \texttt{masses\char`_calc(G)}  returns the vector $[a,b,c,d]$ where:\ps

-- $a$ is the vector of all $[P,m]$ with $P$ in ${\rm Car}_n$ and $m = {\rm m}_{{\rm O}(L)}(P)$ with $m \neq 0$,\ps
-- $b$ is the isomorphism class of the root system ${\rm R}(L)$,\ps
-- $c$ is the vector $(c_j)_{j \in J}$ where $c_j$ encodes both the cycle decomposition of $\gamma_j$ on ${\rm Irr}(R)$ and the integers $d_S$ for each $S$ in ${\rm Irr}_j$,\ps
-- $d$ is the vector $(d_j)_{j \in J}$ with $d_j=[P_j,m_j/M]$ and $M=\sum_{j \in J} m_j$.

\section{The characteristic masses of Niemeier Lattices with roots} \label{par:niemeier}${}^{}$

The aim of this section is to explain a way to determine the characteristic masses of the Niemeier lattices with roots 
which does not use the computationally heavy steps {\sff B1} and {\sff B4} in {\sff Algorithm B}, by rather determining directly the information of step {\sff B6} (and then using of course the elementary results of \S \ref{par:rootlat}). We will use for this the case by case descriptions of these lattices given by Venkov \cite{venkov} or Conway and Sloane \cite[Ch. 16]{conwaysloane}, based on the 
classical connections between lattices and codes \cite{conwaysloane,ebeling}, and study their automorphism groups in slightly more details than what we could find in the literature. To keep this section short, we assume some familiarity with Niemeier lattices and mostly follow the exposition in \cite[Chap. 2.3]{chlannes} to which we refer for more details. \ps\ps

\subsection{Linking modules, Venkov modules and even unimodular lattices}\label{par:linking}${}^{}$
 (a) A {\it  (quadratic)  linking} module\footnote{Such a module is also called a ${\rm qe}$-module in \cite[Chap. 2.3]{chlannes}.} is a finite abelian group $A$ equipped with a quadratic map ${\rm q} : A \longrightarrow \Q/\Z$ whose associated symmetric $\Z$-bilinear map ${\rm b}(x,y):={\rm q}(x+y)-{\rm q}(x)-{\rm q}(y), \,\,\,A \times A \longrightarrow \Q/\Z,$ is a perfect pairing. The isometry group of $A$ is denoted ${\rm O}(A)$.  If $I \subset A$ is a subgroup, we denote by $I^\bot$ the orthogonal of $I$ with respect to ${\rm b}$. We say that $I$ is {\it isotropic} if we have ${\rm q}(I) = 0$ (this is usually stronger than $I \subset I^\bot$). We say that $I$ is a {\it Lagrangian} if it is isotropic and if we have $I=I^\bot$ (or equivalently $|A|=|I|^2$).  \ps\ps
  
 (b) A {\it Venkov module} is a linking module $A$ equipped with a (set theoretic) map ${\rm qm} : A \rightarrow \Q_{\geq 0}$ such that for all $a \in A$ we have ${\rm qm}(a) \equiv {\rm q}(a) \bmod \Z$, ${\rm qm}(0)=0$ and ${\rm qm}(a)>0$ for $a \neq 0$. Venkov modules form an additive category ${\rm Ven}$ in an obvious way; in particular we have an obvious notion of orthogonal direct sum of such objects, denoted $\oplus$. A {\it root} of a Venkov module $A$ is an element $a \in A$ such that ${\rm qm}(a)=1$. \ps\ps
  
(c) Assume $L$ is an even lattice in the Euclidean space $V$. Recall that we set ${\rm q}(x) = \frac{x \cdot x}{2}$ for $x \in V$. The finite abelian group $L^\sharp/L$, equipped with the well-defined quadratic map (that we shall still denote by ${\rm q}$)
$L^\sharp/L \longrightarrow \Q/\Z, \,\, \, x + L \mapsto {\rm q}(x) \bmod \Z$,  is a linking module that we shall denote by ${\rm res}\, L$ (sometimes also called the {\it discriminant group} or {\it glue group} of $L$). This linking module has a canonical structure of Venkov module defined by $${\rm qm}(x) \,=\, {\rm inf}_{y \in x+L} \,\,{\rm q}(y).$$ 
Let $\pi: L^\sharp \rightarrow {\rm res}\, L$ be the canonical projection. The map $I \mapsto \pi^{-1}I$ is a bijection between the set of isotropic subspaces $I$ of ${\rm res}\, L$ and the set of even lattices of $V$ containing $L$. In this bijection, $\pi^{-1}\, I$ is unimodular if, and only if, $I$ is a Lagrangian. Moreover, we have ${\rm R}(\pi^{-1} I)={\rm R}(L)$ if, and only if, $I$ does not contain any root of ${\rm res}\, L$.
\ps\ps

(d) We now focus on the case $L={\rm Q}(R)$ with $R$ an ${\bf ADE}$ root system in $V$. In this case ${\rm Q}(R)^\sharp$ is called the {\it weight} lattice of $R$ and we set ${\rm res}\, R = {\rm res}\,L$. The group ${\rm O}(R)$ naturally acts on ${\rm res}\, R$, with ${\rm W}(R)$ acting trivially, so we have a morphism ${\rm G}(R) \rightarrow {\rm Aut}_{\rm Ven}({\rm res}\, R)$. The Venkov module ${\rm res}\, R$ is the orthogonal sum of the ${\rm res}\, S$ with $S$ an irreducible component of $R$. Assume now $R$ is irreducible. Canonical representatives for the nonzero elements of ${\rm res}\, R$ are given by the so-called {\it minuscule weights} of $R$, that we denote by $\varpi_i$ following Bourbaki's conventions \cite{bourbaki} for the indices. A key property is ${\rm qm}(\varpi_i+{\rm Q}(R))={\rm q}(\varpi_i)$ (see Table \ref{tab:venkov}). We also have ${\rm G}(R) \isomo {\rm Aut}_{\rm Ven}({\rm res}\, R)$. In particular, the element $-{\rm id}$ of ${\rm O}(R)$ is in ${\rm W}(R)$ if, and only if, ${\rm res}\, R$ is a $\Z/2$-vector space.  \ps\ps

\begin{table}[H]
{\scriptsize \renewcommand{\arraystretch}{1.8} \medskip
\begin{center}
\smallskip
\begin{tabu}{c|[1pt]c|c|c|c|c|c}
$R$ &  ${\bf A}_n$ & ${\bf D}_{n}$, $n$\,\,{\rm even} & ${\bf D}_{n}$, $n$\,\,{\rm odd}& ${\bf E}_6$ & ${\bf E}_7$ & ${\bf E}_8$ \cr
\hline
${\rm res}\, R$ & $\Z/(n+1)$ & $\Z/2\times \Z/2$ & $\Z/4$ & $\Z/3$ & $\Z/2$ & $0$ \cr
${\rm min. wts}$ 	& $\varpi_i$, $i=1,\dots,n$ & $\varpi_{n},\varpi_1,\varpi_{n-1}$ & $\varpi_{n},\varpi_1,\varpi_{n-1}$ & $\varpi_1,\varpi_6$ & $\varpi_7$ &  \cr
${\rm class}$ & $i \bmod n+1$ & $\omega, 1, \overline{\omega}$ & $1,2,3 \bmod 4$ & $1,2 \bmod 3$ & $1 \bmod 2$ & \cr
${\rm qm}$ & $\frac{i(n+1-i)}{2(n+1)}$ & $\frac{n}{8},\frac{1}{2},\frac{n}{8}$ &  $\frac{n}{8},\frac{1}{2},\frac{n}{8}$ & $\frac{2}{3}, \frac{2}{3}$ & $\frac{3}{4}$ & \cr
\end{tabu} 
\end{center}
} 
\label{tab:venkov}
\caption{The Venkov module ${\rm res}\, R$ for $R$ an irreducible ${\bf ADE}$ root system.}
\end{table}

\begin{remark} \label{caseD4} {\rm In the case $R \simeq {\bf D}_{n}$ with $n$ even, some authors (e.g. \cite{conwaysloane}) identify the $\Z/2$-vector space ${\rm res}\, R$ with the finite field $\mathbb{F}_4=\{0,1,\omega,\overline{\omega}\}$. Using this identification, the automorphism group of ${\rm res}\, R$, in Ven, is generated by the Frobenius $f(x) = x^2$, as well as $m(x)= \omega x$ for $n=4$ (triality). }
\end{remark}

\subsection{The Niemeier lattices with roots}${}^{}$

Niemeier and Venkov have shown that $L \mapsto {\rm R}(L)$ induces a bijection between the isomorphism classes of Niemeier lattices with roots, and the isomorphism classes of equi-Coxeter\footnote{A root system $R$ is called equi-Coxeter if its irreducible components have the same Coxeter number, then called the Coxeter number of $R$ and denoted ${\rm h}(R)$. The Coxeter numbers of ${\bf A}_n$, ${\bf D}_n$, ${\bf E}_6$, ${\bf E}_7$ and ${\bf E}_8$ are respectively $n+1$, $2n-2$, $12$, $18$ and $30$.} ${\bf ADE}$ root systems in $\R^{24}$. Fix such a root system $R$ in $\R^{24}$. By \S \ref{par:linking}, there is thus a unique ${\rm O}(R)$-orbit of Lagrangians $I$ in ${\rm res}\, R$ containing no root (the "codes"). For any such $I$, the associated Niemeier lattice with root system $R$ is $L = \pi^{-1} I$, we have $|{\rm res}\, R|=|I|^2$ and ${\rm G}(L)={\rm O}(L)/{\rm W}(R)$ is the stabilizer of $I \subset {\rm res}\, R$ in ${\rm G}(R)$. In particular, the conjugacy class of ${\rm G}(L)$ in ${\rm G}(R)$ does not depend on the choice of $I$. \ps\ps

{\bf Goal:} {\it For each of the $23$ possible isomorphism classes of $R$, determine the ${\rm G}(R)$-conjugacy class of the elements of ${\rm G}(L)$ (with their multiplicity)}.\ps\ps

This is exactly the information actually needed to apply Proposition \ref{prop:massnonirr} to each coset $\sigma {\rm W}(R)$ in ${\rm O}(L)$. We will use information on ${\rm G}(L)$ given by Venkov \cite{venkov} and Conway-Sloane \cite[Table 16.1]{conwaysloane} (see also \cite{erokhin2}), such as their order and a composition series. Note that those ${\rm G}(L)$ are also exactly the {\it umbral groups} studied in \cite{umbral}.\footnote{Although we will not use it, as this not the information we need, let us mention that the character tables of umbral groups have been listed in the appendix 2 {\it loc. cit.} (and computed using \texttt{GAP}).}  We may assume that the decomposition of $R$ as a union of its irreducible components has the form\footnote{This is a short notation for $\coprod_{i=1}^g \coprod_{j=1}^{N_i} R_i$.}
$$R \,\,=\,\, N_1R_1 \,\,N_2 R_2\, \, \dots \,\, N_g R_g,$$ 
with $R_i \not \simeq R_j$ for $i \neq j$, and $N_i \geq 1$ for all $i$. We have natural decompositions
$${\rm res}\, R \, \, =\, \, \bigoplus_{i} ({\rm res}\, R_i)^{N_i},\,\, {\rm G}(R) \, \, = \, \, \prod_i \,{\rm G}(N_i R_i) \, \, \text{and} \, \, {\rm G}(N_i R_i) \,=\, {\rm G}(R_i) \,\wr \,{\rm S}_{N_i}.$$
By \eqref{eq:outdyn}, each ${\rm G}(N_i R_i)$ is naturally isomorphic either to the symmetric group ${\rm S}_{N_i}$, to the hyperoctahedral group ${\rm H}_{N_i}$, or to ${\rm T}_{N_i}:={\rm S}_3 \wr {\rm S}_{N_i}$ in the exceptional case $R_i \simeq {\bf D}_4$. The natural exact sequence $1 \rightarrow \prod_i {\rm G}(R_i)^{N_i} \rightarrow {\rm G}(R) \rightarrow \prod_i {\rm S}_{N_i} \rightarrow 1$ induces an exact sequence $1 \rightarrow {\rm G}_1(L) \rightarrow {\rm G}(L) \rightarrow {\rm G}_2(L) \rightarrow 1$.
The orders of ${\rm G}_1(L)$ and ${\rm G}_2(L)$ are given in \cite[Table 16.1]{conwaysloane}. Moreover, the image of ${\rm G}_2(L)$ in ${\rm S}_{N_i}$ is always a transitive subgroup for each $i$.\ps\ps

We denote by $\eta \in {\rm G}(L)$ the class of the element $- {\rm id}$ of ${\rm O}(L)$. It is a central element which does not depend on the choice of $I$, and satisfies $\eta^2=1$. Its image in ${\rm G}(N_i R_i)$ is trivial if $R_i$ has type ${\bf A}_1$, ${\bf D}_{2n}$, ${\bf E}_7$ or ${\bf E}_8$, and equal to the element $-1$ of ${\rm G}(N_i R_i)={\rm H}_{N_i}$ otherwise (\S \ref{par:defhn} and Table \ref{tab:venkov}).  An inspection of Table \cite[Table 16.1]{conwaysloane} shows that we always have ${\rm G}_1(L) = \langle \eta \rangle$, except in the case $R \simeq {\bf D}_4^6$ for which we have ${\rm G}_1(L) \simeq \Z/3$ (and $\eta=1$).  
 \ps\ps

{\sc Notation:} A conjugacy class $C \subset {\rm G}(R)$ has the form $\prod_i C_i$ where $C_i$ is a conjugacy class in $G_i={\rm G}(N_i R_i)$. So $C$ is uniquely determined by the collection $(t_i)$ where $t_i$ is {\it the type} of $C_i$ : a partition of $N_i$ in the case $G_i={\rm S}_{N_i}$, a couple of partitions as in \S \ref{par:conjhn} in the case $G_i={\rm H}_{N_i}$, and similarly a triple of partitions in the case $G_i={\rm T}_{N_i}$. In this last case, and as in \S \ref{par:conjhn}, we use the sequence of symbols 
$ \cdots i^{a_i}\, {\color{cyan} i^{b_i}}\,{\color{magenta} i^{c_i}} \cdots$ 
to denote the conjugacy class whose elements have a cycle decomposition with $a_i$ (resp. $b_i$, $c_i$) cycles of length $i$ whose $i$-th power has order $1$ (resp. $2$, $3$), with same conventions as {\it loc.cit}.\ps\ps

We now start the description of the ${\rm G}(R)$-conjugacy classes of the elements of ${\rm G}(L)$. In the non trivial cases, we list their type and give the number of elements of any given type divided by $|{\rm G}(L)|$ (the {\it size} of the type): \ps\ps
\begin{itemize}
\item $R \simeq {\bf D}_{24},\,\,\, {\bf D}_{16}\,{\bf E}_8,\,\,\, {\bf A}_{24}, \,\,\,{\bf A}_{17}\,{\bf E}_7, \,\,\,{\bf A}_{15}\,{\bf D}_9$ and ${\bf A}_{11}\,{\bf D}_7\,{\bf E}_6$. We have ${\rm G}_2(L)=1$, so ${\rm G}(L)={\rm G}_1(L)= \langle \eta \rangle$. \ps 
\item $R \simeq 3{\bf E}_8$. We have ${\rm G}(R) = {\rm S}_3$ and ${\rm res}\, R = 0$, so ${\rm G}(L)={\rm G}(R)={\rm S}_3$. \ps 
\item $R \simeq 2{\bf D}_{12}$. We have ${\rm G}(R) = {\rm H}_2$, ${\rm G}_1(L)=1$ and ${\rm G}_2(L)={\rm S}_2$. We may take for $I$ the subgroup $\{ 0, (1,\omega), (\omega,1),(\overline{\omega},\overline{\omega})\}$ (note ${\rm q}(I)=\{0, 2,3\}$). For this $I$, ${\rm G}(L)$ is the natural subgroup ${\rm S}_2$ of ${\rm H}_2$. \ps 
\item $R \simeq {\bf D}_{10}\, 2{\bf E}_7$. We have ${\rm G}(R) = {\rm H}_1 \times {\rm H}_2$, ${\rm G}_1(L)=1$ and ${\rm G}_2(L)={\rm S}_2$. We may take $I= \{0, (\omega,1,0),(\overline{\omega},0,1),(1,1,1)\}$ (note ${\rm qm}(I)=\{0,2\}$).
So ${\rm G}(L) \simeq \Z/2$ is generated by the element $( \varepsilon_1 , (1\,2))$ of ${\rm G}(R)$, whose type is (${\color{cyan} 1}$,$2$). \ps 
\item $R \simeq 3{\bf D}_8$. We have ${\rm G}(R) = {\rm H}_3$, ${\rm G}_1(L)=1$ and ${\rm G}_2(L)={\rm S}_3$. We may take 
for $I$ the subgroup generated by the ${\got S}_3$-orbit of $(1,1,\omega)$ (it contains $(0,\overline{\omega},\overline{\omega})$, $(\omega,\omega,\omega)$ and we have $|I|=8$ and ${\rm qm}(I) = \{0,2,3\}$). For this $I$, ${\rm G}(L)$ is the natural subgroup ${\rm S}_3$ of ${\rm H}_3$. \ps 
\item $R \simeq 2 {\bf A}_{12}$. We have ${\rm G}(R) = {\rm H}_2$,  $\eta=-1$ and ${\rm G}_2(L)={\rm S}_2$. We may take $I = \langle a \rangle \simeq \Z/13$ with $a=(1,5)$ (note ${\rm qm}(I) = \{0,2,3\}$). The order $4$ element $\sigma=\varepsilon_1 (1\,2)$ of ${\rm H}_2$ satisfies $\sigma a=-5 a$, hence generates ${\rm G}(L)$. The type of the elements of ${\rm G}(L)$ are thus $1^2$, ${\color{cyan} 1^2}$ and ${\color{cyan} 2}$, with respective size $1/4,1/4$ and $1/2$. \ps 
\item $R \simeq 4{\bf E}_6$. We have ${\rm G}(R) = {\rm H}_4$, $\eta=-1$ and ${\rm G}_2(L)={\rm S}_4$: ${\rm G}(L)$ is a central extension of ${\rm S}_4$ by $\Z/2$. Let $I$ be the Lagrangian of ${\rm res}\, R \simeq (\Z/3)^4$ with $\pi^{-1}I=L$. The stabilizer of $I$ in ${\rm O}({\rm res}\, R)$ is the semi-direct product of ${\rm GL}(I) \simeq {\rm GL}_2(\Z/3)$ and of a $\Z/3$-vector space. The natural morphism ${\rm G}(L) \rightarrow {\rm GL}(I)$ is thus injective, hence bijective. In particular, ${\rm G}(L)$ is not isomorphic to ${\rm S}_4 \times \Z/2$ and we are in case (ii) (c) of Proposition \ref{2Sn}:  ${\rm G}(L)$ is ${\rm H}_4$-conjugate to the subgroup of Example \ref{exgl2z3}. The ${\rm H}_4$-conjugacy classes of ${\rm G}(L)$ are thus given by Table \ref{tablegl2z3}. \ps 
\item $R \simeq 4{\bf D}_6$. We have ${\rm G}(R) = {\rm H}_4$,  ${\rm G}_1(L)=1$ and ${\rm G}_2(L)={\rm S}_4$. We claim that ${\rm G}(L)$ does not contain the natural subgroup ${\rm S}_4$ of ${\rm H}_4$. Indeed, assume that the Lagrangian $I$ of ${\rm res}\, R = ({\rm res}\, {\rm D}_6)^4$ defining $L$ is stable under ${\rm S}_4$. For any $x=(x_1,x_2,x_3,x_4) \in I$ and any $1\leq i \neq j \leq 4$, we have $x+\, (i\,j)x$ in $I$, hence $2 {\rm qm}(x_i+x_j)$ is either $0$ or an integer $\geq 2$. This forces $x_i+x_j=0$ by Table \ref{tab:venkov}, hence $|I| \leq 4$: a contradiction. By Proposition \ref{2Sn}, the subgroup ${\rm G}(L) \subset {\rm H}_4$ is thus ${\rm H}_4$-conjugate to the subgroup ${\rm S}_4^{\epsilon}$, and we are done. \ps 
\item $R \simeq 2{\bf A}_9\,{\bf D}_6$. We have ${\rm G}(R) = {\rm H}_2 \times {\rm H}_1$,  $\eta=(-1,1)$ and ${\rm G}_2(L)={\rm S}_2 \times {\rm S}_1$. We may take for $I$ the subgroup generated by the elements $a=(2,4,0)$, $b=(5,0,\omega)$ and $c=(0,5,\overline{\omega})$, of respective orders $5$, $2$ and $2$ (check ${\rm qm}(I)=\{0,2,3\}$). Observe that $I$ is stable by the element $\sigma = ( \varepsilon_2 (1\,2), -1)$ of ${\rm G}(R)$: we have $\sigma(a) = (4,-2,0)=2a$, $\sigma(b)=(0,-5,\overline{\omega})=c$ and $\sigma(c)=(5,0,\omega)=b$. This shows ${\rm G}(L) = \langle \sigma \rangle \simeq \Z/4$, the types of its elements being $(1^2,1)$, $({\color{cyan} 1^2},1)$ and $({\color{cyan} 2},{\color{cyan} 1})$, with respective size $1/4, 1/4$ and $1/2$. \ps 
\item $R \simeq 3{\bf A}_8$. We have ${\rm G}(R) = {\rm H}_3$, $\eta=-1$ and ${\rm G}_2(L)={\rm S}_3$.  By Proposition \ref{2Sn}, ${\rm G}(L)$ is ${\rm H}_3$-conjugate to the subgroup $\{\pm 1\} \cdot {\rm S}_3$ of ${\rm H}_3$. \ps 
\item $R \simeq 2{\bf A}_7\,2{\bf D}_5$. We have ${\rm G}(R) = {\rm H}_2 \times {\rm H}_2$, $\eta=(-1,-1)$ and ${\rm G}_2(L)={\rm S}_2 \times {\rm S}_2$. We may take for $I$ the subgroup generated by 
the elements $a=(1,1,1,2)$ and $b=(1,-1,2,1)$ of order $8$ (check ${\rm qm}(I)=\{0,2,3\}$). Note that $I$ is stable under $\sigma_1=((1\, 2),\varepsilon_2)$ and $\sigma_2=(\varepsilon_2,(1\,2))$: we have $\sigma_1(a)=a$, $\sigma_1(b)=(-1,1,2,3)=-b$, and $\sigma_2$ exchanges $a$ and $b$. This shows ${\rm G}(L)=\langle \sigma_1,\sigma_2 \rangle$ (dihedral of order $8$), with types $(1^2,1^2)$, $({\color{cyan} 1^2},{\color{cyan}1^2})$ of size $1/8$, and types $(2,1\,{\color{cyan}1})$, $(1\,{\color{cyan}1},2)$ and $({\color{cyan}2},{\color{cyan}2})$ of size $1/4$.
\ps 
\item $R \simeq 4{\bf A}_6$. We have ${\rm G}(R) = {\rm H}_4$,  $\eta=-1$ and $|{\rm G}_2(L)|=12$. So ${\rm G}(L)$ is a central extension of ${\rm G}_2(L)={\rm Alt}_{4}$ by $\Z/2$. It has an injective morphism to ${\rm GL}(I)={\rm GL}_2(\Z/7)$ (same argument as for $4{\bf E}_6$): this is a non split extension. By Remark \ref{rem:A4}, the types of the elements of ${\rm G}(L)$ follow thus from Table \ref{tablegl2z3}: they are $1^4$, ${\color{cyan} 1^4}$, ${\color{cyan} 2^2} $, $1\,3$ and ${\color{cyan}1\,3 }$, with respective sizes $1/24$, $1/24$, $1/4$, $1/3$ and $1/3$. 
\ps 
\item $R \simeq 6{\bf D}_4$. We have ${\rm G}(R) = {\rm T}_6$,  $|{\rm G}_1(L)|=3$ and ${\rm G}_2(L)={\rm S}_6$. We identify ${\rm res}\, {\rm D}_4$ with $\mathbb{F}_4$ as in Remark~\ref{caseD4}. Following Conway and Sloane, $I$ is an {\it hexacode} in $\mathbb{F}_4^6$. By \cite[\S 11.2]{conwaysloane}, we may choose for $I$ the $\mathbb{F}_4$-vector space generated by the $K$-orbit of 
$(\omega,\overline{\omega},\omega,\overline{\omega},\omega,\overline{\omega})$, where 
$K$ is the subgroup of ${\rm Alt}_6$ preserving $\{\{1,2\},\{3,4\},\{5,6\}\}$. For this choice, ${\rm G}(L)$ contains $K$ and ${\rm G}_1(L)$ is generated by the element $(m,m,m,m,m,m)$ of ${\rm Aut}_{\Z/2}(\mathbb{F}_4)^6=({\rm S}_3)^6$. Two other elements of ${\rm G}(L)$ are for instance $(f,f,f,f,f,f) (12)$ and $(1,1,1,1,m^2,m)(1\,2\,3)$ (for the latter, recall that $(\omega,\omega,\overline{\omega},\overline{\omega},1,1)$ is in $I$). As $K$, $(1\,2)$ and $(1\,2\,3)$ generate ${\rm S}_6$, a straightforward computation allows to list the types of the elements of ${\rm G}(L)$: we obtain Table~\ref{table3S6}.
\begin{table}[H]
{\scriptsize \renewcommand{\arraystretch}{1.8} \medskip
\begin{center}
\begin{tabu}{c|c|c|c|c|c|c|c|c}
${\texttt{type}}$ & $1^6$ & ${\color{magenta} 1^6}$ & ${\color{cyan} 1^4}\,2$ & $1^2\,2^2$ & ${\color{magenta}1^2\,2^2}$ & $2^3$ & $1\,{\color{magenta} 1^2}\,3$ & $3^2$ \\
\hline
${\texttt{size}}$  & $1/2160$ & $1/1080$ & $1/48$ & $1/48$ & $1/24$ & $1/48$ & $1/18$ & $1/18$ \\
\noalign{\hrule height 1pt}
${\texttt{type}}$ & ${\color{cyan}1^2}\, 4$ & $2\,4$ & ${\color{magenta} 2\,4}$ & ${\color{magenta} 1\,5}$ & $1\,5$ & $6$ & ${\color{cyan}1}\,{\color{magenta}2}\,{\color{cyan}3}$  & \\
\hline
${\texttt{size}}$  & $1/8$  & $1/24$  & $1/12$ & $2/15$ & $1/15$ & $1/6$ & $1/6$ &
\end{tabu} 
\end{center}} 
\caption{{\small The ${\rm T}_6$-conjugacy classes of the elements of $3. {\rm S}_6$.} }\label{table3S6}
\end{table}
\item $R \simeq 4{\bf A}_5\,{\bf D}_4$. We have ${\rm G}(R) = {\rm H}_4 \times {\rm T}_1$,  $\eta=(-1,1)$ and ${\rm G}_2(L)={\rm S}_4 \times {\rm S}_1$. We identify ${\rm res}\, {\rm D}_4$ with $\mathbb{F}_4$ as in Remark \ref{caseD4}. The first projection ${\rm pr}_1: {\rm G}(L) \rightarrow {\rm H}_4$ is injective, and its image ${\rm H}(L)$ is a central extension of ${\rm S}_4$ by $\Z/2$. By Proposition~\ref{2Sn}, ${\rm H}(L)$ is either conjugate to $\{\pm 1\} \times {\rm S}_4$ or to the group ${\rm GL}_2(\Z/3)$ embeded as in Example~\ref{exgl2z3}. By \cite[\S 11.2]{conwaysloane}, we may take for $I$ the subgroup generated by the $\sigma$-orbit of $a=(2,0,2,4,0)$ and $b=(3,3,0,0,\overline{w})$, where $\sigma=((2\, 3\, 4), m)$.  In order to determine ${\rm H}(L)$ it is enough to find the unique $\pm v \in \{\pm 1\}^4$ such that $v (1\,2) \in {\rm H}(L)$. Note that $2I$ is the $\Z/3$-vector space generated by $(2,0,2,4,0)$ and $(2,4,0,2,0)$. We deduce $v = \pm (1,1,1,-1)$: we have ${\rm G}(L)\simeq {\rm H}(L) \simeq {\rm GL}_2(\Z/3)$. On the other hand, the second projection ${\rm pr}_2: {\rm G}(L) \rightarrow {\rm T}_1$ is trivial on $\eta$ hence factors through a morphism $\mu : {\rm S}_4 \rightarrow {\rm S}_3$. We have $\sigma \in {\rm G}(L)$ and ${\rm pr}_2(\sigma) = m$ has order $3$: $\mu$ is ``the'' classical surjective morphism from ${\rm S}_4$ to ${\rm S}_3$. The types of ${\rm G}(L)$ are thus immediately deduced from Table~\ref{tablegl2z3}. \ps 
\item $R \simeq 6{\bf A}_4$. We have ${\rm G}(R) = {\rm H}_6$,  $\eta=-1$ and ${\rm G}_2(L)$ is a transitive subgroup of ${\rm S}_6$ of order $120$, so ${\rm G}_2(L)$ is isomorphic to ${\rm S}_5$ and ${\rm G}(L)$ is a central extension of ${\rm S}_5$ by $\Z/2$ in ${\rm H}_6$. We claim that ${\rm G}(L)$ does not contain the triple transposition $\tau = (12)\,(34)\,(56)$. Indeed, otherwise the Lagrangian $I$ defining $L$ would be invariant by $\tau$. Set $I^{\pm}=\{x \in I, \tau(x) = \pm x\}$.
A nonzero element of $I^+$ has the form $(a,a,b,b,c,c)$ with $2{\rm qm}(a)+2{\rm qm}(b)+2{\rm qm}(c)$ an integer $\neq 1$. This forces $\{\pm a,\pm b,\pm c\}=\Z/5$ since ${\rm qm}({\rm res}\, {\rm A}_4) = \{ 0, \frac{2}{5},\frac{3}{5}\}$. But the nondegenerate conic $a^2+b^2+c^2 =0$ in $(\Z/5)^3$ contains all those vectors: we have $\dim_{\Z/5} I^+ \leq 1$. A similar argument shows $\dim_{\Z/5} I^- \leq 1$. This is a contradiction as $I = I^+ \oplus I^-$ has dimension $3$, hence the claim. By Proposition \ref{corS5} and the discussion after this proposition, ${\rm G}(L)$ is ${\rm H}_6$-conjugate to the image of the map \eqref{iotaS5}. The type of its elements are thus given by Table \ref{tab:2S5}. \ps 
\item $R \simeq 8{\bf A}_3$. We have ${\rm G}(R) = {\rm H}_8$,  $\eta=-1$ and $|{\rm G}_2(L)|=1344$. 
In this case, $I \subset {\rm res}\, R = (\Z/4)^8$ is a so-called {\it octacode} \cite{conwaysloane}. 
The subgroup $C:=I/2I$ a Hamming code in ${\rm res}\, R \otimes \Z/2 = (\Z/2)^8$ and ${\rm G}_2(L)$ is the automorphism group of this code. In particular, $C$ is included in the hyperplane $H$ of $(\Z/2)^8$ defined by $\sum_i x_i = 0$ and $V:=H/C$ is a hyperplane in the $4$-dimensional $\Z/2$-vector space $W:=(\Z/2)^8/C$. The (easy) theory of Hamming codes shows that the map $\iota : \{1,\dots,8\} \rightarrow W,$ sending $j$ to the class of the canonical basis element $\delta_j$ of $(\Z/2)^8$, is injective with image an affine hyperplane under $V$, and identifies ${\rm G}_2(L)$ with the affine group of $\{1,\dots,8\}$ for this affine structure. In particular, ${\rm G}_2(L)$ is isomorphic to ${\rm GA}_3(\Z/2)= (\Z/2)^3 \rtimes {\rm GL}_3(\Z/2)$. To go further we choose some $I$: following \cite[Table 16.1]{conwaysloane} we take the subgroup generated by the $c$-orbit of the element $(3,2,0,0,1,0,1,1)$ where $c$ is the $7$-cycle $(2\,3\,4\,5\,6\,7\,8)$ (we have ${\rm qm}(I)=\{0,2,3,4\}$). With this choice of $I$, we have $c \in {\rm G}(L)$ and checks that $\tau=(3\,4\,6)(5\,8\,7)$ and $\sigma=\varepsilon_3 \varepsilon_6  \varepsilon_7 \varepsilon_8\,(2\,3) (4\,5\,6\,8)$ lie in ${\rm G}(L)$ as well. The images in ${\rm G}_2(L)$ of $1,c,c^{-1},\tau,\sigma$ and $\sigma^2$ are representatives of the conjugacy classes of the stabilizer ${\rm G}_2(L)_1$ of $1$ in $\{1,\dots,8\}$, with resp. sizes $1/168$, $1/7$, $1/7$, $1/3$, $1/4$ and $1/8$ (recall ${\rm G}_2(L)_1 \simeq {\rm GL}_3(\Z/2)$). But $c,\tau$ and $\sigma$ belong to the stabilizer ${\rm G}(L)_1$ of $1 \in \{\pm 1\}^8$ in ${\rm G}(L)$: the natural map ${\rm G}(L) \rightarrow {\rm G}_2(L)$ induces an isomorphism ${\rm G}(L)_1 \isomo {\rm G}_2(L)_1$. An inspection of $C$ shows that the translation by the class of $\delta_1-\delta_2 $ in $V$ is the element $(1\,2)\,(3\,7)\,(4\,5)\,(6\,8)$ of ${\rm S}_8$. One deduces from these information representatives of the conjugacy classes of ${\rm G}(L)$: their types are gathered in Table \ref{table2ga3z2}.

\begin{table}[H]
{\scriptsize \renewcommand{\arraystretch}{1.8} \medskip
\begin{center}
\begin{tabular}{c|c|c|c|c|c|c|c|c|c|c|c|c}
${\texttt{type}}$ &  {\color{cyan} $1^8$} & $1^2${\color{cyan}$1^2$}$2^2$ & $2^4$ & {\color{cyan}$2^4$}& $1^2 3^2$ &{\color{cyan}$1^2 3^2$} &$1\,${\color{cyan}$1\,2\,$}$4$& $4^2$ & {\color{cyan}$4^2$} & $2\,6\,$& $1\,7\,$ & {\color{cyan}$1\,7\,$}\\
\hline   \texttt{size} & $1/2688$ & $1/32$ & $1/192$ & $1/32$ & $1/12$ & $1/12$ & $1/8$ & $1/16$ & $1/8$ & $1/6$ & $1/7$ & $1/7$
\end{tabular} 
\end{center}} 
\caption{{\small The ${\rm H}_8$-conjugacy classes of the nontrivial elements of $2. {\rm GA}_3(\Z/2)$.} }\label{table2ga3z2}
\end{table}

\item $R \simeq 12{\bf A}_2$. We have ${\rm G}(R) = {\rm H}_{12}$,  $\eta=-1$ and ${\rm G}_2(L)$ is isomorphic to the Mathieu group ${\rm M}_{12}$. The Lagrangian $I \subset {\rm res}\,R\, = (\Z/3)^{12}$ is a {\it ternary Golay code}, whose automorphism group ${\rm G}(L)$ is the central extension of ${\rm M}_{12}$ by $\Z/2$ denoted $2.{\rm M}_{12}$ in the $\mathbb{ATLAS}$.  We know since Frobenius \cite[p. 11]{frobenius} the cycle decompositions, and cardinality, of all the conjugacy classes of ${\rm M}_{12}$. The inverse image in $2.{\rm M}_{12}$ of such a class $c$ is the union of one or two conjugacy classes $c' \cup -c'$, the cycle decomposition of $c'$ being the same as that of $c$ except that each cycle of $c$ now has a sign to be determined. It is an amusing exercise\footnote{That such an exercise is possible follows from the following fact: if we have an equality of polynomials
$\prod_i (t^i-1)^{a_i}(t^i+1)^{b_i}=\prod_i (t^i-1)^{a'_i}(t^i+1)^{b'_i}$ with $a_i+b_i=a'_i+b'_i$ for each $i$, then $a_i=a'_i$ and $b_i=b'_i$ for each $i$ (use $t^i+1=(t^{2i}-1)/(t^i-1)$ and Lemma~\ref{uniquenesscyc}).} to extract these signs from the lines $\chi_2$ and $\chi_{18}$, and from the power maps, of the character table of $2.{\rm M}_{12}$ in the $\mathbb{ATLAS}$. We obtain Table~\ref{table2M12}.\footnote{An alternative way to proceed is to use the description of $2. {\rm M}_{12}$ given by Hall in \cite{hallM12}, as the automorphism group of a $12 \times 12$ Hadamard matrix (a subgroup of ${\rm H}_{12}$). Using the $4$ generators given by Hall {\it loc. cit.}, and applying \texttt{GAP}'s \texttt{ConjugacyClasses} algorithm to the permutation group on $24$ letters they generate, we confirm Table~\ref{table2M12}.}

\begin{table}[H]
{\scriptsize \renewcommand{\arraystretch}{1.8} \medskip
\begin{center}
\begin{tabular}{c|c|c|c|c|c|c|c|c|c|c}
${\texttt{type}}$ & {\color{cyan}$1^{12}$} & {\color{cyan}$2^6$} & $1^4 2^4$ & {\color{cyan}$1^{4}$}$2^4$ & $1^3 3^3$ & {\color{cyan}$1^3 3^3$} & $3^4$ & {\color{cyan}$3^4$} & $2^2 4^2$ & $1^2${\color{cyan}$1^2$}$4^2$ \\
\hline   \texttt{size} & $1/190080$ & $1/240$ & $1/384$ & $1/384$ & $1/108$ & $1/108$ & $1/72$ & $1/72$ & $1/32$ & $1/32$ \\
\noalign{\hrule height 1pt}
${\texttt{type}}$ & $1^2 5^2$ & {\color{cyan}$1^2 5^2$}  & {\color{cyan}$6^2$}  & $1\,2\,3\,6$ & {\color{cyan}$1\,$}$2\,${\color{cyan}$3\,$}$6$ & $4\,8$ & $1\,${\color{cyan}$1\,2\,$}$8$ & {\color{cyan}$2\,10$}& $1\,11$ & {\color{cyan}$1\,11$}  \\ 
\hline \texttt{size} & $1/20$ & $1/20$  & $1/12$ & $1/12$ & $1/12$ & $1/8$ & $1/8$ & $1/10$ & $1/11$ & $1/11$  
\end{tabular} 
\end{center}} 
\caption{{\small The ${\rm H}_{12}$-conjugacy classes of the nontrivial elements of $2. {\rm M}_{12}$.} }\label{table2M12}
\end{table}

\item $R \simeq 24{\bf A}_1$. We have ${\rm G}(R) = {\rm S}_{24}$ and ${\rm G}(L)$ is a Mathieu group ${\rm M}_{24}$. The cycle decompositions and cardinality of the conjugacy classes of ${\rm M}_{24}$ are given by Frobenius in \cite[p. 12-13]{frobenius}: see Table \ref{tableM24}. \ps 

\begin{table}[htp]
{\scriptsize \renewcommand{\arraystretch}{1.8} \medskip
\begin{center}
\begin{tabular}{c|c|c|c|c|c|c|c|c|c|c}
${\texttt{type}}$ & $1^8\,2^8$ & $2^{12}$ & $1^6\,3^6$ & $3^8$ & $2^4\,4^4$  & $1^4\,2^2\,4^4$ & $4^6$ & $1^4\,5^4$ & $1^2\,2^2\,3^2\,6^2$ & $6^4$    \\
\hline
${\texttt{mass}}$ & $1/21504$ & $1/7680$ & $1/1080$ & $1/504$ & $1/384$ & $1/128$ & $1/96$ & $1/60$ & $1/24$  & $1/24$ \\
\noalign{\hrule height 1pt}
${\texttt{type}}$& $1^3 7^3$ & $1^2\,2\,4\,8^2$ & $2^2\,10^2$  & $1^2\,11^2$ & $2\,4\,6\,12$ & $12^2$ & $1\,2\,7\,14$ & $1\,3\,5\,15$ & $3\,21$ & $1\,23$ \\
\hline
${\texttt{mass}}$  & $1/21$ & $1/16$ & $1/20$ & $1/11$ & $1/12$ & $1/12$ & $1/7$ & $2/15$ & $2/21$ & $2/23$ \\
\end{tabular}
\end{center}}
\caption{{\small The ${\rm S}_{24}$-conjugacy classes of  the nontrivial elements of ${\rm M}_{24}$.} }\label{tableM24}
\end{table}

\end{itemize}

{\bf Comparison with the output of {\sff Algorithm B}}. For each of the $23$ root systems $R$ above, we verified that the types and sizes of the ${\rm G}(R)$-conjugacy classes of ${\rm G}(L)$ found are exactly those returned (from scratch, and in a few seconds!) by {\sff Algorithm B} (components $3$ and $4$ returned by \texttt{masses\char`_calc}, see \S \ref{par:implementalgo}). The natural isomorphism ${\rm O}(L)_\rho \simeq {\rm G}(L)$ and {\sff Algorithm B} provide thus a rather useful tool to study the groups ${\rm G}(L)$.

\clearpage

\begin{appendix}

\titleformat{\section}{\bfseries}{\appendixname~\thesection .}{0.5em}{}
\titleformat{\subsection}{\bfseries}{\thesection .~\arabic{subsection}}{0.5em}{}
\renewcommand{\thesubsection}{\Alph{section}.~\arabic{subsection}}

\section{Irreducible characters of compact orthogonal groups}\label{par:on}${}^{}$

Let $n\geq 1$ be an integer. We denote by ${\rm O}(n)$ the isometry group of the standard Euclidean space $V=\R^n$. 
We know since Weyl that the complex, irreducible, continuous representations of the compact group ${\rm O}(n)$ 
are all defined over $\R$ and parameterized in a natural way by the {\it $n$-permissible (integer) partitions} $\lambda$. In this section, we recall this parameterization and discuss
formulas for the irreducible characters due to Weyl and Koike-Terrada. \ps

\subsection{The $n$-permissible partitions} ${}^{}$

Recall that a partition $\lambda$ is a non-increasing integer sequence  $\lambda_1 \geq \lambda_2 \geq \cdots$ 
with $\lambda_i \geq 0$ for all $i\geq 1$ and $\lambda_i =0$ for $i$ big enough. 
We also say that $\lambda$ is a partition of the integer $|\lambda| := \sum_i \lambda_i$. The {\it diagram} of $\lambda$ is the Young diagram whose $i$-th row has $\lambda_i$ 
boxes for each $i \geq 1$. The {\it dual} of $\lambda$ is the partition $\lambda^\ast$ defined by $\lambda^\ast_i = |\{j \geq 1 \,|\,\lambda_j \geq i\}|$ (with ``transpose'' diagram). \ps

Following Weyl, the partition $\lambda$ is called {\it $n$-permissible} if the first two columns of its diagram contain at most $n$ boxes, or equivalently if we have $\lambda^\ast_1+ \lambda^\ast_2 \leq n$. 
If $\lambda$ is $n$-permissible, there is a unique $n$-permissible partition $\mu$ with $\lambda^\ast_i=\mu^\ast_i$ for $i>1$ and $\lambda^\ast_1+\mu^\ast_1=n$, called the {\it associate} of $\lambda$ and denoted ${\rm ass}(\lambda)$. The map $\lambda \mapsto {\rm ass}(\lambda)$ is an involution of the set of $n$-permissible integer partitions. \ps
An partition $\lambda$ is called {\it $n$-positive} if we have
$\lambda^\ast_1 \leq n/2$ (hence $\lambda_i=0$ for $i>n/2$). If $\lambda$ is $n$-admissible but not $n$-positive, then ${\rm ass}(\lambda)$ is $n$-positive. \ps

\subsection{Weyl's construction}\label{par:weylconst}${}^{}$

 For any integer $d\geq 0$, we consider following Weyl the kernel ${\rm K}_d(V)$ of the direct sum of the $d(d-1)/2$ contraction maps\footnote{All tensor products are taken over $\R$ in \S \ref{par:weylconst}.} ${\rm c}_{i,j} : V^{\otimes d} \rightarrow V^{\otimes (d-2)}$, defined for $1 \leq i<j \leq d$ by ${\rm c}_{i,j}(v_1 \otimes v_2 \otimes \cdots \otimes v_d)\,=\,(v_i \cdot v_j) \,v_1 \otimes v_2 \otimes \cdots \otimes \widehat{v_i} \otimes \cdots \otimes \widehat{v_j} \otimes \cdots \otimes v_d$. This kernel has a natural linear action of ${\rm O}(n) \times \got{S}_d$, hence decomposes as
$${\rm K}_d(V) \simeq \bigoplus_{\{\lambda\,\,| \, \, |\lambda|=d\}}\,  {\rm K}_\lambda(V) \otimes {\rm R}_\lambda$$
where ${\rm R}_\lambda$ is ``the'' irreducible representation of $\got{S}_d$ classically parameterized by $\lambda$, and ${\rm K}_\lambda(V)$ is a real representation of ${\rm O}(n)$.  Set ${\rm W}_\lambda = {\rm K}_\lambda(V) \otimes \C$. \ps

Weyl shows that ${\rm W}_\lambda$ is either $0$ or an irreducible representation of ${\rm O}(n)$  \cite[Thm. 5.7.D]{weyl}. Moreover, ${\rm W}_\lambda$  is nonzero if and only if $\lambda$ is $n$-permissible \cite[Thm. 5.7.A \& C]{weyl}. Moreover, he shows that $\lambda \rightarrow {\rm W}_\lambda$ is a bijection between the set of $n$-permissible partitions and the isomorphism classes of irreducible representations of ${\rm O}(n)$ \cite[Thm. 5.7.H \& 7.9.B]{weyl}. 
The element $-{\rm 1}_n$ clearly acts as multiplication by $(-1)^d$ on ${\rm W}_\lambda$. Weyl shows 
\begin{equation} \label{twistdet} {\rm W}_{{\rm ass}(\lambda)} \simeq {\rm W}_\lambda \otimes \det\end{equation}
and studies the restriction of ${\rm W}_\lambda$  to the index two subgroup ${\rm SO}(n) \subset {\rm O}(n)$ in Chap. V.9 \& VII.9.
We may assume $\lambda$ is $n$-positive. There are two cases:\ps

(i) $\lambda \neq {\rm ass}(\lambda)$. The restriction of ${\rm W}_\lambda$ to ${\rm SO}(n)$ is then irreducible 
with highest weight $\sum_{i \leq n/2} \,\,\lambda_i \,\varepsilon_i$, using the classical notations of Bourbaki \cite[Pl. IV]{bourbaki}. Moreover, the natural action 
of ${\rm O}(n)/{\rm SO}(n)=\Z/2$ on the highest weight lines of ${\rm W}_\lambda$ is trivial (and non trivial on those of ${\rm W}_{{\rm ass}(\lambda)} \otimes \C$). \ps

(ii) $\lambda = {\rm ass}(\lambda)$. This forces $n \equiv 0 \bmod 2$ and $\lambda_{n/2}>0$. The restriction of ${\rm W}_\lambda$ to ${\rm SO}(n)$ is then 
the sum of the two irreducible representations, conjugate under ${\rm O}(n)$, with highest weights $\pm \lambda_{n/2} \,\varepsilon_{n/2}\,+\,\sum_{i=1}^{n/2-1} \lambda_i \,\varepsilon_i$. \ps

\subsection{Character formulas} ${}^{}$ 

Weyl gives a determinantal formula for the character of ${\rm W}_\lambda$ in \cite[Theorem 7.9.A]{weyl}.
Contrary to the standard so-called {\it Weyl character formula}, which applies to any connected compact Lie groups, that formula equally applies
to elements in any of the two connected components\footnote{Let us mention that there exists also a variant of the Weyl character formula which applies to the irreducible characters of non connected 
compact Lie groups: see e.g. \cite{kostant,wendt}.} of ${\rm O}(n)$. Assume $g$ is in ${\rm O}(n)$ and write $\det(1-tg)^{-1} = \sum_{i \in \Z} p_i t^i$ in $\Z[[t]]$  (so $p_i=0$ for $i<0$). 
Weyl shows {\it loc. cit.} that for any $n$-permissible partitions $\lambda$ we have
\begin{equation}\label{weylcarf} {\rm Trace}(g \, ;\, {\rm W}_\lambda) = \det \,(p_{{\lambda_i}-i+j}- p_{\lambda_i-i-j})_{1 \leq i,j \leq \lambda_1^\ast}.\end{equation}
If we write $\det(1+tg)=  \sum_{i \in \Z} e_i t^i$ (so $e_i=0$ for $i<0$ or $i>n$), and set $\delta_1=0$ and $\delta_j=1$ for $j>1$, then \cite[Theorem 2.3.3 (6)]{koiketerada} implies 
\begin{equation}\label{univid} \det \,(p_{{\lambda_i}-i+j}- p_{\lambda_i-i-j})_{1 \leq i,j \leq \lambda_1^\ast} =  \det (e_{\lambda_i^\ast-i+j}\,+\, \delta_{j}\,e_{\lambda_i^\ast-i-j+2})_{1 \leq i,j \leq \lambda_1},\end{equation}
See also the equivalence of (ii) and (iv) in \cite[Cor. A.46]{fh} for a direct alternative proof of this equality. \ps

\begin{remark} {\rm In the case $\lambda_1$=0, or equivalently $|\lambda|=0$ or $\lambda$ is the empty diagram, then ${\rm W}_{\lambda}$ is the trivial representation and both determinants above are indeed $1$ by convention. 
Moreover, the formula $e_{n-i} = (\det g) e_i$ for $i \in \Z$ shows that the determinant on the right-hand side of \eqref{univid} is multiplied by $\det g$ if $\lambda$ is replaced by ${\rm ass}(\lambda)$ (it amounts to multiply by $\det g$ the first line of the matrix inside the determinant), in agreement with Formula \eqref{twistdet}.
}
\end{remark}

\section{An asymptotic formula}\label{par:asymptotic} 

\begin{prop}\label{prop:asymp} Let $L$ be a lattice in the Euclidean space $\R^n$ and $\lambda$ an $n$-permissible partition with $|\lambda| \equiv 0 \bmod 2$. Then we have $$\dim {\rm W}_\lambda^{{\rm O}(L)} \sim \frac{2 }{|{\rm O}(L)|} \dim {\rm W}_\lambda$$
for $\lambda \rightarrow \infty$, in the sense that $\lambda_i - \lambda_{i+1} \rightarrow +\infty$ for each $1 \leq i \leq n/2$. 
\end{prop}

\begin{pf} As we have $\lambda \rightarrow \infty$ we may assume $\lambda$ is positive and $\lambda_{[n/2]} >0$.
Denote by ${\rm V}_\lambda$ the irreducible constituent of $({\rm W}_\lambda)_{|{\rm SO}(n)}$ with highest weight $\sum_{i \leq n/2} \,\,\lambda_i \,\varepsilon_i$. Set ${\rm SO}(L)={\rm O}(L) \cap {\rm SO}(n)$. If $n$ is odd, we have $({\rm W}_\lambda)_{|{\rm SO}(n)} = {\rm V}_\lambda$, ${\rm O}(L) = \{ \pm {\rm id}\} \times {\rm SO}(L)$ and ${\rm W}_\lambda^{{\rm O}(L)}={\rm V}_\lambda^{{\rm SO}(L)}$. If $n$ is even, then $({\rm W}_\lambda)_{|{\rm SO}(n)}$ is the direct sum of ${\rm V}_\lambda$ and of its outer conjugate ${\rm V}_{\lambda}'$, and ${\rm W}_\lambda$ is induced from ${\rm V}_\lambda$: we have thus ${\rm W}_\lambda^{{\rm O}(L)} = {\rm V}_\lambda^{{\rm SO}(L)}$ in the case ${\rm O}(L) \neq {\rm SO}(L)$ and ${\rm W}_\lambda^{{\rm O}(L)} = {\rm V}_\lambda^{{\rm SO}(L)} \oplus ({\rm V}_\lambda')^{{\rm SO}(L)}$ otherwise.
We conclude from the degenerate form of Weyl's character formula for ${\rm SO}(n)$ given in \cite[Prop. 1.9]{chcl}.\end{pf}

Assume now $n \equiv -1,0,1 \mod 8$ and set $\mu_n = \sum_{[L] \in {\rm X}_n} \frac{1}{|{\rm O}(L)|}$. The mass formula of Minkowski-Siegel-Smith asserts that we have $\mu_n = |\frac{{\rm B}_{n/2}}{n}\prod_{j=1}^{n/2-1}\frac{{\rm B}_{2j}}{4j}|$ for $n \equiv 0 \bmod 8$, and $\mu_n =|\prod_{j=1}^{(n-1)/2}\frac{{\rm B}_{2j}}{4j}|$ for $n \equiv \pm 1 \bmod 8$, where the ${\rm B}_{m}$ are the Bernouilli numbers  \cite{conwaysloanemass}. \ps\ps

\begin{cor} For $n \equiv -1,0,1 \bmod 8$, $|\lambda| \equiv 0 \bmod 2$ and $\lambda \rightarrow \infty$ we have
$\dim {\rm M}_{{\rm W}_\lambda}({\rm O}_n) \,\,\sim \,\,2\, \mu_n \,\dim {\rm W}_\lambda$.
\end{cor} 
\ps\ps
For instance, in the case $n=24$ of main interest here we have $\mu_{24}  \approx 8 \cdot 10^{-15}$, quite a small number compared to $|{\rm X}_{24}|=24$, and of course we expect $\dim {\rm M}_{W_\lambda} ({\rm O}_{24})$ to be small for small values of $\lambda$.  	

\newpage

\section{Tables}\label{sec:tables}

\begin{table}[H]
{\scriptsize
\renewcommand{\arraystretch}{1.8} \medskip
\begin{center}
\scalebox{.9}{
\begin{tabu}{c|c|[1pt]c|c|[1pt]c|c|[1pt]c|c|[1pt]c|c|[1pt]c|c|[1pt]c|c}$P$ & $m$ & $P$ & $m$ & $P$ & $m$ & $P$ & $m$ & $P$ & $m$ & $P$ & $m$ & $P$ & $m$ \\
\hline
$3^{2}$ & $1/24$ & 
$6^{2}$ & $1/24$ & 
$1^{2}3$ & $1/12$ & 
$2^{2}6$ & $1/12$ & 
$1\,2\,3$ & $1/4$ & 
$12$ & $1/4$ & 
$1\,2\,6$ & $1/4$ \\
\end{tabu}}\end{center}}
\caption{{\footnotesize The $7$ nonzero $m={\rm m}_{\sigma {\rm W}({\rm D}_4)}(P)$ for $P$ in ${\rm Car}_{4}$, where $\sigma$ in ${\rm G}({\rm D}_4)$ has order $3$.}}\label{TableCWD4}
\end{table}

\begin{table}[H]
{\tiny
\renewcommand{\arraystretch}{1.8} \medskip
\scalebox{.9}{
\hspace{-1.7cm}\begin{tabu}{c|c|[1pt]c|c|[1pt]c|c|[1pt]c|c|[1pt]c|c|[1pt]c|c|[1pt]c|c|[1pt]c|c|[1pt]c|c}$P$ & $m$ & $P$ & $m$ & $P$ & $m$ & $P$ & $m$ & $P$ & $m$ & $P$ & $m$ & $P$ & $m$ & $P$ & $m$ & $P$ & $m$ \\
\hline
$1^{6}$ & $1/51840$ & 
$3^{3}$ & $1/648$ & 
$1^{2}3^{2}$ & $1/108$ & 
$1\,2^{3}4$ & $1/96$ & 
$1\,2\,3^{2}$ & $1/36$ & 
$1^{3}2\,4$ & $1/32$ & 
$3\,12$ & $1/12$ & 
$1^{2}5$ & $1/10$ & 
$1\,2\,8$ & $1/8$ \\
$1^{5}2$ & $1/1440$ & 
$1^{4}3$ & $1/216$ & 
$1^{2}4^{2}$ & $1/96$ & 
$3\,6^{2}$ & $1/72$ & 
$1^{2}2^{2}6$ & $1/36$ & 
$1^{2}2^{2}3$ & $1/24$ & 
$1\,2\,3\,6$ & $1/12$ & 
$1\,2\,5$ & $1/10$ & 
 & $$ \\
$1^{2}2^{4}$ & $1/1152$ & 
$1^{4}2^{2}$ & $1/192$ & 
$1^{3}2^{3}$ & $1/96$ & 
$1^{3}2\,3$ & $1/36$ & 
$2^{2}3\,6$ & $1/36$ & 
$1^{2}2^{2}4$ & $1/16$ & 
$1\,2\,4\,6$ & $1/12$ & 
$9$ & $1/9$ & 
 & $$ \\
\end{tabu}
}
\caption{{\footnotesize The $25$ nonzero $m={\rm m}_{{\rm W}({\rm E}_6)}(P)$ for $P$ in ${\rm Car}_6$.}}\label{TableWE6}}
\end{table}

\vspace{-.5cm}

\begin{table}[H]
{\tiny
\renewcommand{\arraystretch}{1.8} \medskip
\scalebox{.9}{
\hspace{-2.5cm}
\begin{tabu}{c|c|[1pt]c|c|[1pt]c|c|[1pt]c|c|[1pt]c|c|[1pt]c|c|[1pt]c|c|[1pt]c|c|[1pt]c|c}
$P$ & $m$ & $P$ & $m$ & $P$ & $m$ & $P$ & $m$ & $P$ & $m$ & $P$ & $m$ & $P$ & $m$ & $P$ & $m$ & $P$ & $m$ \\
\hline
$1^{7}$ & $1/2903040$ & 
$1^{5}2^{2}$ & $1/3072$ & 
$1^{4}2^{3}$ & $13/9216$ & 
$1^{3}2^{2}6$ & $1/288$ & 
$2\,3^{2}6$ & $1/144$ & 
$1^{3}5$ & $1/60$ & 
$1^{2}2\,4\,6$ & $1/48$ & 
$1\,3\,5$ & $1/30$ & 
$1\,9$ & $1/18$ \\
$2^{7}$ & $1/2903040$ & 
$1^{2}2^{5}$ & $1/3072$ & 
$1^{3}2^{4}$ & $13/9216$ & 
$1\,2^{4}6$ & $1/288$ & 
$1\,3\,6^{2}$ & $1/144$ & 
$2^{3}10$ & $1/60$ & 
$1\,2^{2}4\,6$ & $1/48$ & 
$2\,6\,10$ & $1/30$ & 
$2\,18$ & $1/18$ \\
$1^{6}2$ & $1/46080$ & 
$1\,3^{3}$ & $1/1296$ & 
$1^{4}2\,4$ & $1/384$ & 
$1^{2}2\,4^{2}$ & $1/256$ & 
$1^{3}2^{2}3$ & $1/96$ & 
$1^{3}2^{2}4$ & $7/384$ & 
$1^{2}2\,8$ & $1/32$ & 
$1\,3\,12$ & $1/24$ & 
$1^{2}2\,3\,6$ & $1/18$ \\
$1\,2^{6}$ & $1/46080$ & 
$2\,6^{3}$ & $1/1296$ & 
$1\,2^{4}4$ & $1/384$ & 
$1\,2^{2}4^{2}$ & $1/256$ & 
$1^{2}2^{3}6$ & $1/96$ & 
$1^{2}2^{3}4$ & $7/384$ & 
$1\,4\,8$ & $1/32$ & 
$2\,6\,12$ & $1/24$ & 
$1\,2^{2}3\,6$ & $1/18$ \\
$1^{5}3$ & $1/4320$ & 
$1^{3}4^{2}$ & $1/768$ & 
$1^{4}2\,3$ & $1/288$ & 
$1^{3}3^{2}$ & $1/216$ & 
$1^{2}2\,3^{2}$ & $1/72$ & 
$1^{2}2\,3\,4$ & $1/48$ & 
$1\,2^{2}8$ & $1/32$ & 
$1^{2}2\,5$ & $1/20$ & 
$1\,7$ & $1/14$ \\
$2^{5}6$ & $1/4320$ & 
$2^{3}4^{2}$ & $1/768$ & 
$1^{2}2^{3}3$ & $1/288$ & 
$2^{3}6^{2}$ & $1/216$ & 
$1\,2^{2}6^{2}$ & $1/72$ & 
$1\,2^{2}3\,4$ & $1/48$ & 
$2\,4\,8$ & $1/32$ & 
$1\,2^{2}10$ & $1/20$ & 
$2\,14$ & $1/14$ \\
\end{tabu}}
\caption{{\footnotesize The $54$ nonzero $m={\rm m}_{{\rm W}({\rm E}_7)}(P)$ for $P$ in ${\rm Car}_7$.}}\label{TableWE7}}
\end{table}
\vspace{-.7cm}
\begin{table}[H]
{\tiny
\renewcommand{\arraystretch}{1.8} \medskip
\scalebox{.9}{
\hspace{-2cm}
\begin{tabu}{c|c|[1pt]c|c|[1pt]c|c|[1pt]c|c|[1pt]c|c|[1pt]c|c|[1pt]c|c|[1pt]c|c}
$P$ & $m$ & $P$ & $m$ & $P$ & $m$ & $P$ & $m$ & $P$ & $m$ & $P$ & $m$ & $P$ & $m$ & $P$ & $m$ \\
\hline
$1^{8}$ & $1/696729600$ & 
$2^{4}4^{2}$ & $1/18432$ & 
$3^{2}6^{2}$ & $1/1728$ & 
$1^{2}3\,4^{2}$ & $1/576$ & 
$1\,2\,3\,6^{2}$ & $1/288$ & 
$1\,2^{3}10$ & $1/120$ & 
$2^{2}6\,10$ & $1/60$ & 
$1^{2}7$ & $1/28$ \\
$2^{8}$ & $1/696729600$ & 
$1^{5}2\,4$ & $1/15360$ & 
$1^{3}2\,4^{2}$ & $1/1536$ & 
$1^{3}2^{3}3$ & $1/576$ & 
$1^{2}2^{2}6^{2}$ & $1/288$ & 
$1^{2}2^{2}4^{2}$ & $9/1024$ & 
$3\,9$ & $1/54$ & 
$1\,2\,7$ & $1/28$ \\
$1^{7}2$ & $1/5806080$ & 
$1\,2^{5}4$ & $1/15360$ & 
$1\,2^{3}4^{2}$ & $1/1536$ & 
$1\,2^{3}3\,4$ & $1/576$ & 
$6^{2}12$ & $1/288$ & 
$1^{2}9$ & $1/108$ & 
$6\,18$ & $1/54$ & 
$1\,2\,14$ & $1/28$ \\
$1\,2^{7}$ & $1/5806080$ & 
$1^{5}2\,3$ & $1/8640$ & 
$1^{4}5$ & $1/1200$ & 
$1^{3}2\,4\,6$ & $1/576$ & 
$1^{3}2^{3}4$ & $19/4608$ & 
$2^{2}18$ & $1/108$ & 
$1\,2\,3\,8$ & $1/48$ & 
$2^{2}14$ & $1/28$ \\
$1^{6}3$ & $1/311040$ & 
$1\,2^{5}6$ & $1/8640$ & 
$2^{4}10$ & $1/1200$ & 
$1^{3}2^{3}6$ & $1/576$ & 
$8^{2}$ & $1/192$ & 
$1^{2}2^{2}3\,4$ & $1/96$ & 
$1\,2\,3\,12$ & $1/48$ & 
$24$ & $1/24$ \\
$2^{6}6$ & $1/311040$ & 
$1^{2}3^{3}$ & $1/7776$ & 
$1^{4}2^{2}3$ & $1/1152$ & 
$2^{2}4^{2}6$ & $1/576$ & 
$1^{3}2\,3\,4$ & $1/192$ & 
$1^{2}2^{2}4\,6$ & $1/96$ & 
$1\,2\,6\,8$ & $1/48$ & 
$1\,2\,4\,12$ & $1/24$ \\
$1^{6}2^{2}$ & $1/184320$ & 
$2^{2}6^{3}$ & $1/7776$ & 
$1^{2}2^{4}6$ & $1/1152$ & 
$1^{3}2\,3^{2}$ & $1/432$ & 
$1\,2^{3}4\,6$ & $1/192$ & 
$1^{2}2^{2}5$ & $1/80$ & 
$1\,2\,6\,12$ & $1/48$ & 
$20$ & $1/20$ \\
$1^{2}2^{6}$ & $1/184320$ & 
$1^{2}2^{4}3$ & $1/6912$ & 
$2^{2}3^{2}6$ & $1/864$ & 
$1\,2^{3}6^{2}$ & $1/432$ & 
$1^{2}3\,12$ & $1/144$ & 
$1^{2}2^{2}10$ & $1/80$ & 
$1\,2\,4\,5$ & $1/40$ & 
$1\,2\,4\,8$ & $5/64$ \\
$3^{4}$ & $1/155520$ & 
$1^{4}2^{2}6$ & $1/6912$ & 
$1^{2}3\,6^{2}$ & $1/864$ & 
$1^{3}2\,8$ & $1/384$ & 
$1^{3}2\,3\,6$ & $1/144$ & 
$4^{2}12$ & $1/72$ & 
$1\,2\,4\,10$ & $1/40$ & 
 & $$ \\
$6^{4}$ & $1/155520$ & 
$1^{4}2^{4}$ & $37/221184$ & 
$1^{4}2^{2}4$ & $1/768$ & 
$1\,2^{3}8$ & $1/384$ & 
$1\,2^{3}3\,6$ & $1/144$ & 
$1\,2\,3\,4\,6$ & $1/72$ & 
$1\,2\,9$ & $1/36$ & 
 & $$ \\
$4^{4}$ & $1/46080$ & 
$1^{4}3^{2}$ & $1/2592$ & 
$1\,2\,4^{3}$ & $1/768$ & 
$1^{2}2^{2}3^{2}$ & $1/288$ & 
$2^{2}6\,12$ & $1/144$ & 
$1^{2}2^{2}8$ & $1/64$ & 
$1\,2\,18$ & $1/36$ & 
 & $$ \\
$1^{4}4^{2}$ & $1/18432$ & 
$1\,2\,3^{3}$ & $1/2592$ & 
$1^{2}2^{4}4$ & $1/768$ & 
$3^{2}12$ & $1/288$ & 
$1^{2}4\,8$ & $1/128$ & 
$1^{2}3\,5$ & $1/60$ & 
$1^{2}2^{2}3\,6$ & $1/36$ & 
 & $$ \\
$1^{5}2^{3}$ & $1/18432$ & 
$1\,2\,6^{3}$ & $1/2592$ & 
$5^{2}$ & $1/600$ & 
$12^{2}$ & $1/288$ & 
$2^{2}4\,8$ & $1/128$ & 
$1\,2\,3\,5$ & $1/60$ & 
$15$ & $1/30$ & 
 & $$ \\
$1^{3}2^{5}$ & $1/18432$ & 
$2^{4}6^{2}$ & $1/2592$ & 
$10^{2}$ & $1/600$ & 
$1\,2\,3^{2}6$ & $1/288$ & 
$1^{3}2\,5$ & $1/120$ & 
$1\,2\,6\,10$ & $1/60$ & 
$30$ & $1/30$ & 
 & $$ \\
\end{tabu}
}
\caption{{\footnotesize The $106$ nonzero $m={\rm m}_{{\rm W}({\rm E}_8)}(P)$ for $P$ in ${\rm Car}_8$.}}\label{TableWE8}
} 
\end{table}

\newpage

\begin{table}[H]
\label{TableLeech}
{\tiny
\renewcommand{\arraystretch}{1.8} \medskip
\caption{{\footnotesize The $160$ nonzero $m={\rm m}_{{\rm O}({\rm Leech})}(P)$ for $P$ in ${\rm Car}_{24}$.}}
\medskip
\scalebox{.9}{
\hspace{-2.5cm}
\begin{tabu}{c|c|[1pt]c|c|[1pt]c|c|[1pt]c|c|[1pt]c|c|[1pt]c|c}
$P$ & $m$ & $P$ & $m$ & $P$ & $m$ & $P$ & $m$ & $P$ & $m$ & $P$ & $m$ \\
\hline
$1^{24}$ & $1/8315553613086720000$ & 
$2^{6}3^{4}6^{5}$ & $1/311040$ & 
$3^{4}6^{4}12^{2}$ & $1/4608$ & 
$1^{4}2^{4}4^{4}8^{2}$ & $1/768$ & 
$4^{2}8^{2}12\,24$ & $1/144$ & 
$1^{2}3\,5\,15^{2}$ & $1/60$ \\
$2^{24}$ & $1/8315553613086720000$ & 
$1^{8}2^{8}4^{4}$ & $1/294912$ & 
$3^{3}9^{3}$ & $1/3888$ & 
$12^{2}24^{2}$ & $1/576$ & 
$1^{4}11^{2}$ & $1/132$ & 
$1^{2}3\,5\,15\,30$ & $1/60$ \\
$3^{12}$ & $1/2690072985600$ & 
$3^{4}12^{4}$ & $1/276480$ & 
$6^{3}18^{3}$ & $1/3888$ & 
$1^{4}2^{4}3^{3}4^{2}6\,12$ & $1/576$ & 
$2^{4}22^{2}$ & $1/132$ & 
$2^{2}6\,10\,15\,30$ & $1/60$ \\
$6^{12}$ & $1/2690072985600$ & 
$6^{4}12^{4}$ & $1/276480$ & 
$1^{4}3^{2}4^{4}12^{2}$ & $1/3456$ & 
$1^{4}2^{2}3^{3}6^{2}12^{2}$ & $1/576$ & 
$1^{2}4\,8\,16^{2}$ & $1/128$ & 
$1^{2}2^{2}3\,5\,6\,10\,15$ & $1/60$ \\
$1^{16}2^{8}$ & $1/178362777600$ & 
$4^{4}8^{4}$ & $1/92160$ & 
$1^{6}2^{6}3^{3}6^{3}$ & $1/3456$ & 
$2^{2}3^{2}4^{2}6\,12^{3}$ & $1/576$ & 
$2^{2}4\,8\,16^{2}$ & $1/128$ & 
$2^{2}6\,10\,30^{2}$ & $1/60$ \\
$1^{8}2^{16}$ & $1/178362777600$ & 
$1^{8}5^{4}$ & $1/72000$ & 
$2^{4}4^{4}6^{2}12^{2}$ & $1/3456$ & 
$1^{2}2^{2}3\,4^{4}6\,12^{2}$ & $1/576$ & 
$7^{2}21$ & $1/126$ & 
$1^{2}2^{2}3\,5\,6\,10\,30$ & $1/60$ \\
$4^{12}$ & $1/2012774400$ & 
$2^{8}10^{4}$ & $1/72000$ & 
$1^{6}7^{3}$ & $1/2352$ & 
$1^{2}3\,4^{2}6^{2}12^{3}$ & $1/576$ & 
$14^{2}42$ & $1/126$ & 
$1^{2}4^{2}7\,28$ & $1/56$ \\
$1^{8}4^{8}$ & $1/743178240$ & 
$8^{6}$ & $1/48384$ & 
$2^{6}14^{3}$ & $1/2352$ & 
$1^{4}2^{4}3^{4}6^{4}$ & $1/576$ & 
$15\,60$ & $1/120$ & 
$2^{2}4^{2}14\,28$ & $1/56$ \\
$2^{8}4^{8}$ & $1/743178240$ & 
$1^{8}2^{4}3^{4}6^{2}$ & $1/41472$ & 
$1^{2}2^{2}4^{2}8^{4}$ & $1/2048$ & 
$1^{2}2^{4}3^{2}6^{3}12^{2}$ & $1/576$ & 
$30\,60$ & $1/120$ & 
$52$ & $1/52$ \\
$1^{12}2^{12}$ & $1/389283840$ & 
$1^{4}2^{8}3^{2}6^{4}$ & $1/41472$ & 
$3^{4}15^{2}$ & $1/1800$ & 
$1^{4}2^{4}3\,4^{2}6^{3}12$ & $1/576$ & 
$1^{4}2^{2}7^{2}14$ & $1/112$ & 
$1^{2}2^{2}3\,6\,8^{2}24$ & $1/48$ \\
$1^{12}3^{6}$ & $1/117573120$ & 
$7^{4}$ & $1/35280$ & 
$6^{4}30^{2}$ & $1/1800$ & 
$21^{2}$ & $1/504$ & 
$1^{2}2^{4}7\,14^{2}$ & $1/112$ & 
$1^{2}2^{2}3^{2}4\,6\,12\,24$ & $1/48$ \\
$2^{12}6^{6}$ & $1/117573120$ & 
$14^{4}$ & $1/35280$ & 
$1^{6}2^{4}4^{3}8^{2}$ & $1/1536$ & 
$42^{2}$ & $1/504$ & 
$1^{2}2^{2}3\,9^{2}18$ & $1/108$ & 
$1^{2}2^{2}3\,4\,6^{2}12\,24$ & $1/48$ \\
$1^{6}3^{9}$ & $1/25194240$ & 
$1^{4}2^{2}3^{6}6^{3}$ & $1/31104$ & 
$1^{4}2^{6}4^{3}8^{2}$ & $1/1536$ & 
$3^{2}6\,9^{2}18$ & $1/432$ & 
$1^{2}3\,6\,9^{2}18$ & $1/108$ & 
$1^{2}2^{2}3^{2}4^{2}6^{2}12^{2}$ & $1/48$ \\
$2^{6}6^{9}$ & $1/25194240$ & 
$1^{2}2^{4}3^{3}6^{6}$ & $1/31104$ & 
$4^{4}12^{4}$ & $1/1440$ & 
$3\,6^{2}9\,18^{2}$ & $1/432$ & 
$1^{2}2^{2}6\,9\,18^{2}$ & $1/108$ & 
$1^{2}3^{2}7\,21$ & $1/42$ \\
$3^{8}6^{4}$ & $1/19906560$ & 
$1^{4}5^{5}$ & $1/30000$ & 
$20^{3}$ & $1/1200$ & 
$3^{2}12\,24^{2}$ & $1/384$ & 
$2^{2}3\,6\,9\,18^{2}$ & $1/108$ & 
$2^{2}6^{2}14\,42$ & $1/42$ \\
$3^{4}6^{8}$ & $1/19906560$ & 
$2^{4}10^{5}$ & $1/30000$ & 
$4^{4}20^{2}$ & $1/1200$ & 
$6^{2}12\,24^{2}$ & $1/384$ & 
$3^{2}6^{2}12^{2}24$ & $1/96$ & 
$20\,40$ & $1/40$ \\
$5^{6}$ & $1/6048000$ & 
$5^{4}10^{2}$ & $1/19200$ & 
$1^{4}5^{3}10^{2}$ & $1/1200$ & 
$5^{2}15^{2}$ & $1/360$ & 
$84$ & $1/84$ & 
$1^{2}2^{2}5^{2}10\,20$ & $1/40$ \\
$10^{6}$ & $1/6048000$ & 
$5^{2}10^{4}$ & $1/19200$ & 
$1^{4}2^{4}5^{3}10$ & $1/1200$ & 
$10^{2}30^{2}$ & $1/360$ & 
$3\,9\,12\,36$ & $1/72$ & 
$1^{2}2^{2}5\,10^{2}20$ & $1/40$ \\
$1^{10}2^{6}4^{4}$ & $1/1474560$ & 
$1^{6}2^{6}4^{6}$ & $1/15360$ & 
$2^{4}5^{2}10^{3}$ & $1/1200$ & 
$13^{2}$ & $1/312$ & 
$8^{2}24^{2}$ & $1/72$ & 
$1^{2}2^{2}4^{2}5\,10\,20$ & $1/40$ \\
$1^{6}2^{10}4^{4}$ & $1/1474560$ & 
$1^{4}4^{2}8^{4}$ & $1/12288$ & 
$1^{4}2^{4}5\,10^{3}$ & $1/1200$ & 
$26^{2}$ & $1/312$ & 
$6\,12\,18\,36$ & $1/72$ & 
$39$ & $1/39$ \\
$1^{4}2^{4}4^{8}$ & $1/1179648$ & 
$2^{4}4^{2}8^{4}$ & $1/12288$ & 
$1^{4}3\,9^{3}$ & $1/972$ & 
$1^{4}2^{4}3^{2}4^{2}6^{2}12$ & $1/288$ & 
$35$ & $1/70$ & 
$78$ & $1/39$ \\
$1^{8}3^{8}$ & $1/1088640$ & 
$15^{3}$ & $1/10800$ & 
$1^{2}3^{2}9^{3}$ & $1/972$ & 
$15^{2}30$ & $1/240$ & 
$70$ & $1/70$ & 
$56$ & $1/28$ \\
$2^{8}6^{8}$ & $1/1088640$ & 
$30^{3}$ & $1/10800$ & 
$2^{4}6\,18^{3}$ & $1/972$ & 
$15\,30^{2}$ & $1/240$ & 
$3^{2}33$ & $1/66$ & 
$1^{2}23$ & $1/23$ \\
$12^{6}$ & $1/483840$ & 
$3^{2}6^{2}12^{4}$ & $1/9216$ & 
$2^{2}6^{2}18^{3}$ & $1/972$ & 
$1^{4}3^{2}5^{2}15$ & $1/180$ & 
$6^{2}66$ & $1/66$ & 
$2^{2}46$ & $1/23$ \\
$1^{6}2^{6}3^{5}6$ & $1/311040$ & 
$1^{4}2^{4}8^{4}$ & $1/6144$ & 
$5^{2}20^{2}$ & $1/960$ & 
$2^{4}6^{2}10^{2}30$ & $1/180$ & 
$1^{2}2^{2}8\,16^{2}$ & $1/64$ & 
$1^{2}2^{2}11\,22$ & $1/22$ \\
$1^{6}3^{5}6^{4}$ & $1/311040$ & 
$1^{2}3^{3}4^{2}12^{3}$ & $1/5184$ & 
$10^{2}20^{2}$ & $1/960$ & 
$28^{2}$ & $1/168$ & 
$20\,60$ & $1/60$ & 
 & $$ \\
$1^{6}2^{6}3\,6^{5}$ & $1/311040$ & 
$2^{2}4^{2}6^{3}12^{3}$ & $1/5184$ & 
$24^{3}$ & $1/864$ & 
$1^{4}2^{4}5^{2}10^{2}$ & $1/160$ & 
$12^{2}60$ & $1/60$ & 
 & $$ \\
\end{tabu}
}
}
\end{table}

\newpage

\begin{table}[H]
{\scriptsize \renewcommand{\arraystretch}{1.8} \medskip
\begin{center}
\caption{{\small The nonzero dimensions of  ${\rm M}_{{\rm W}_\lambda}({\rm O}_{24})$ for $\lambda_1 \leq 3$.} }
\medskip
\begin{tabu}{c|c|[1pt]c|c|[1pt]c|c|[1pt]c|c|[1pt]c|c}
$\lambda$ & $\dim $  & $\lambda$ & $\dim$ & $\lambda$ & $\dim $  & $\lambda$ & $\dim $ & $\lambda$ & $\dim $  \cr
\hline
$ $ & $24 : 1$ &
$ 3^{2} $ & $4 : 0$ &
$ 3^{3} 2^{4} 1 $ & $181 : 0$ &
$ 3^{4} 2^{8} $ & $148$ &
$ 3^{6} 2^{6} $ & $276$ \cr \hline
$ 1^{8} $ & $1 : 1$ &
$ 3^{2} 1^{2} $ & $19 : 0$ &
$ 3^{3} 2^{4} 1^{3} $ & $97 : 0$ &
$ 3^{5} 1 $ & $27 : 0$ &
$ 3^{7} 1 $ & $174 : 0$ \cr \hline
$ 1^{12} $ & $1$ &
$ 3^{2} 2 $ & $3 : 0$ &
$ 3^{3} 2^{4} 1^{5} $ & $1$ &
$ 3^{5} 1^{3} $ & $94 : 1$ &
$ 3^{7} 1^{3} $ & $333 : 17$ \cr \hline
$ 2 $ & $9 : 0$ &
$ 3^{2} 2 1^{2} $ & $19 : 0$ &
$ 3^{3} 2^{5} 1 $ & $251 : 1$ &
$ 3^{5} 1^{5} $ & $20 : 0$ &
$ 3^{7} 1^{5} $ & $211$ \cr \hline
$ 2^{2} $ & $27 : 0$ &
$ 3^{2} 2^{2} $ & $15 : 0$ &
$ 3^{3} 2^{5} 1^{3} $ & $120 : 1$ &
$ 3^{5} 1^{7} $ & $1$ &
$ 3^{7} 2 1 $ & $512 : 17$ \cr \hline
$ 2^{3} $ & $26 : 0$ &
$ 3^{2} 2^{2} 1^{2} $ & $50 : 0$ &
$ 3^{3} 2^{6} 1 $ & $265 : 1$ &
$ 3^{5} 2 1 $ & $140 : 0$ &
$ 3^{7} 2 1^{3} $ & $801 : 342$ \cr \hline
$ 2^{4} $ & $43 : 0$ &
$ 3^{2} 2^{3} $ & $18 : 0$ &
$ 3^{3} 2^{6} 1^{3} $ & $100$ &
$ 3^{5} 2 1^{3} $ & $242 : 1$ &
$ 3^{7} 2^{2} 1 $ & $905 : 253$ \cr \hline
$ 2^{4} 1^{8} $ & $1$ &
$ 3^{2} 2^{3} 1^{2} $ & $58 : 0$ &
$ 3^{3} 2^{7} 1 $ & $219 : 51$ &
$ 3^{5} 2 1^{5} $ & $82 : 1$ &
$ 3^{7} 2^{2} 1^{3} $ & $927$ \cr \hline
$ 2^{5} $ & $35 : 0$ &
$ 3^{2} 2^{4} $ & $46 : 0$ &
$ 3^{3} 2^{8} 1 $ & $134$ &
$ 3^{5} 2^{2} 1 $ & $308 : 0$ &
$ 3^{7} 2^{3} 1 $ & $1042 : 683$ \cr \hline
$ 2^{5} 1^{4} $ & $1 : 1$ &
$ 3^{2} 2^{4} 1^{2} $ & $97 : 0$ &
$ 3^{4} $ & $28 : 0$ &
$ 3^{5} 2^{2} 1^{3} $ & $417 : 1$ &
$ 3^{7} 2^{4} 1 $ & $675$ \cr \hline
$ 2^{6} $ & $67 : 1$ &
$ 3^{2} 2^{4} 1^{6} $ & $1$ &
$ 3^{4} 1^{2} $ & $28 : 0$ &
$ 3^{5} 2^{2} 1^{5} $ & $87$ &
$ 3^{8} $ & $191 : 34$ \cr \hline
$ 2^{6} 1^{4} $ & $1 : 1$ &
$ 3^{2} 2^{5} $ & $48 : 0$ &
$ 3^{4} 1^{4} $ & $53 : 1$ &
$ 3^{5} 2^{3} 1 $ & $546 : 1$ &
$ 3^{8} 1^{2} $ & $476 : 137$ \cr \hline
$ 2^{6} 1^{6} $ & $1$ &
$ 3^{2} 2^{5} 1^{2} $ & $91 : 0$ &
$ 3^{4} 1^{8} $ & $1$ &
$ 3^{5} 2^{3} 1^{3} $ & $551 : 111$ &
$ 3^{8} 1^{4} $ & $530$ \cr \hline
$ 2^{7} $ & $42 : 0$ &
$ 3^{2} 2^{6} $ & $97 : 0$ &
$ 3^{4} 2 $ & $30 : 0$ &
$ 3^{5} 2^{4} 1 $ & $672 : 58$ &
$ 3^{8} 2 $ & $327 : 51$ \cr \hline
$ 2^{8} $ & $69 : 1$ &
$ 3^{2} 2^{6} 1^{2} $ & $123 : 1$ &
$ 3^{4} 2 1^{2} $ & $80 : 0$ &
$ 3^{5} 2^{4} 1^{3} $ & $525$ &
$ 3^{8} 2 1^{2} $ & $881 : 552$ \cr \hline
$ 2^{8} 1^{4} $ & $1$ &
$ 3^{2} 2^{6} 1^{4} $ & $1$ &
$ 3^{4} 2 1^{4} $ & $51 : 0$ &
$ 3^{5} 2^{5} 1 $ & $659 : 325$ &
$ 3^{8} 2^{2} $ & $660 : 333$ \cr \hline
$ 2^{9} $ & $37 : 0$ &
$ 3^{2} 2^{7} $ & $70 : 0$ &
$ 3^{4} 2 1^{6} $ & $1 : 1$ &
$ 3^{5} 2^{6} 1 $ & $398$ &
$ 3^{8} 2^{2} 1^{2} $ & $1047$ \cr \hline
$ 2^{10} $ & $48 : 0$ &
$ 3^{2} 2^{7} 1^{2} $ & $74 : 0$ &
$ 3^{4} 2^{2} $ & $112 : 1$ &
$ 3^{6} $ & $36 : 0$ &
$ 3^{8} 2^{3} $ & $500 : 364$ \cr \hline
$ 2^{11} $ & $11 : 0$ &
$ 3^{2} 2^{8} $ & $104 : 0$ &
$ 3^{4} 2^{2} 1^{2} $ & $202 : 1$ &
$ 3^{6} 1^{2} $ & $217 : 1$ &
$ 3^{8} 2^{4} $ & $346$ \cr \hline
$ 2^{12} $ & $37$ &
$ 3^{2} 2^{8} 1^{2} $ & $86$ &
$ 3^{4} 2^{2} 1^{4} $ & $132 : 2$ &
$ 3^{6} 1^{4} $ & $180 : 0$ &
$ 3^{9} 1 $ & $307 : 133$ \cr \hline
$ 3 1 $ & $1 : 0$ &
$ 3^{2} 2^{9} $ & $39 : 8$ &
$ 3^{4} 2^{2} 1^{6} $ & $1$ &
$ 3^{6} 1^{6} $ & $91$ &
$ 3^{9} 1^{3} $ & $496$ \cr \hline
$ 3 2 1 $ & $7 : 0$ &
$ 3^{2} 2^{10} $ & $54$ &
$ 3^{4} 2^{3} $ & $155 : 0$ &
$ 3^{6} 2 $ & $79 : 0$ &
$ 3^{9} 2 1 $ & $651 : 491$ \cr \hline
$ 3 2^{2} 1 $ & $11 : 0$ &
$ 3^{3} 1 $ & $8 : 0$ &
$ 3^{4} 2^{3} 1^{2} $ & $291 : 0$ &
$ 3^{6} 2 1^{2} $ & $474 : 0$ &
$ 3^{9} 2^{2} 1 $ & $542$ \cr \hline
$ 3 2^{3} 1 $ & $31 : 0$ &
$ 3^{3} 1^{3} $ & $6 : 0$ &
$ 3^{4} 2^{3} 1^{4} $ & $126 : 1$ &
$ 3^{6} 2 1^{4} $ & $367 : 61$ &
$ 3^{10} $ & $158 : 121$ \cr \hline
$ 3 2^{4} 1 $ & $33 : 0$ &
$ 3^{3} 2 1 $ & $25 : 0$ &
$ 3^{4} 2^{4} $ & $293 : 1$ &
$ 3^{6} 2^{2} $ & $270 : 0$ &
$ 3^{10} 1^{2} $ & $406$ \cr \hline
$ 3 2^{5} 1 $ & $56 : 0$ &
$ 3^{3} 2 1^{3} $ & $33 : 0$ &
$ 3^{4} 2^{4} 1^{2} $ & $432 : 1$ &
$ 3^{6} 2^{2} 1^{2} $ & $902 : 93$ &
$ 3^{10} 2 $ & $177 : 160$ \cr \hline
$ 3 2^{6} 1 $ & $61 : 0$ &
$ 3^{3} 2^{2} 1 $ & $67 : 0$ &
$ 3^{4} 2^{4} 1^{4} $ & $156$ &
$ 3^{6} 2^{2} 1^{4} $ & $551$ &
$ 3^{10} 2^{2} $ & $161$ \cr \hline
$ 3 2^{7} 1 $ & $63 : 0$ &
$ 3^{3} 2^{2} 1^{3} $ & $49 : 0$ &
$ 3^{4} 2^{5} $ & $270 : 0$ &
$ 3^{6} 2^{3} $ & $386 : 16$ &
$ 3^{11} 1 $ & $93$ \cr \hline
$ 3 2^{8} 1 $ & $59 : 0$ &
$ 3^{3} 2^{2} 1^{5} $ & $1 : 1$ &
$ 3^{4} 2^{5} 1^{2} $ & $387 : 75$ &
$ 3^{6} 2^{3} 1^{2} $ & $988 : 418$ &
$ 3^{12} $ & $74$ \cr \hline
$ 3 2^{8} 1^{3} $ & $1$ &
$ 3^{3} 2^{3} 1 $ & $122 : 0$ &
$ 3^{4} 2^{6} $ & $380 : 73$ &
$ 3^{6} 2^{4} $ & $563 : 197$ &
& \cr \hline
$ 3 2^{9} 1 $ & $53 : 0$ &
$ 3^{3} 2^{3} 1^{3} $ & $102 : 1$ &
$ 3^{4} 2^{6} 1^{2} $ & $362$ &
$ 3^{6} 2^{4} 1^{2} $ & $948$ &
& \cr \hline
$ 3 2^{10} 1 $ & $18$ &
$ 3^{3} 2^{3} 1^{5} $ & $1 : 1$ &
$ 3^{4} 2^{7} $ & $192 : 89$ &
$ 3^{6} 2^{5} $ & $371 : 286$ &
& \cr 
\end{tabu} 
\end{center}
} 
\label{tab:dimO243}
\end{table}

\end{appendix}

\end{document}